\documentclass[preprint,3p,fleqn,sort]{elsarticle}
\usepackage[labelfont=bf, singlelinecheck=false,
            justification=justified]{caption}
\usepackage{xcolor}
\usepackage{amsmath,amssymb}
\usepackage{mathrsfs}
\usepackage{stmaryrd} % Brackets for v_jump
\usepackage{booktabs}
\usepackage{microtype}
\usepackage{mathtools} % shortintertext
\usepackage{multirow}
\usepackage{subfig}
\usepackage{enumitem}
\usepackage{setspace}
\usepackage{placeins} % FloatBarrier
\usepackage{etoolbox,siunitx}
\usepackage{bm}
\usepackage{tikz}
\usepackage[numbers]{natbib}

\robustify\bfseries
\makeatletter
\g@addto@macro\normalsize{%
	\setlength\abovedisplayskip{4pt}
	\setlength\belowdisplayskip{4pt}
	\setlength\abovedisplayshortskip{4pt}
	\setlength\belowdisplayshortskip{4pt}
}
\makeatother
\usepackage[T1]{fontenc}
\usepackage{lmodern}
\usepackage{accsupp} % for ensuring the right Unicode codepoint upon pasting
\SetSymbolFont{stmry}{bold}{U}{stmry}{m}{n}
\usepackage{lineno,hyperref}
\renewcommand{\vec}[1]{\ensuremath{\boldsymbol #1}}
\renewcommand{\vec}[1]{\ensuremath{\boldsymbol #1}}

\newcommand{\halb}{\frac{1}{2}}
\newcommand{\Hlim}{$H^{(c)}_{3L}$}

\newcommand{\sgn}{\ensuremath{\text{sgn}}}
\begin{document}

\begin{frontmatter}

\title{Third-order Limiting for Hyperbolic Conservation Laws applied to Adaptive Mesh Refinement and Non-Uniform 2D Grids}

\author[aachen]{Birte Schmidtmann\corref{mycorrespondingauthor}}
\cortext[mycorrespondingauthor]{Corresponding author}
\ead{schmidtmann@mathcces.rwth-aachen.de}
\author[duesseldorf]{Pawel Buchm\"uller}
\author[aachen]{Manuel Torrilhon}
\address[aachen]{MathCCES, RWTH Aachen, Schinkelstr. 2, 52062 Aachen}
\address[duesseldorf]{Institute of Mathematics, Heinrich-Heine University D\"usseldorf, Universit\"atsstr. 1, 40225 D\"usseldorf, Germany}
%\address[koeln]{Mathematisches Institut, Universit\"at zu K\"oln, Weyertal 86-90, 50931 K\"oln}

\numberwithin{equation}{section}


\begin{abstract}
	In this paper we extend the recently developed third-order limiter function $H_{3\text{L}}^{(c)}$ 
	[\textit{J. Sci. Comput.}, (2016), 68(2), pp.~624--652] 
	%\cite{SchmidtmannSeiboldTorrilhon2015}. 
	to make it applicable for more elaborate test cases in the context of finite volume schemes.
	This work covers the generalization to non-uniform grids in one and two space dimensions, as well as 
	two-dimensional Cartesian grids with adaptive mesh refinement (AMR). The extension to 2D is obtained 
	by the common approach of dimensional splitting. 
	In order to apply this technique without loss of third-order accuracy, the order-fix developed by 
	Buchm\"uller and Helzel [\textit{J. Sci. Comput.}, (2014), 61(2), pp.~343--368] 
	%\cite{BuchmuellerHelzel2014} 
	is incorporated into the scheme. Several numerical examples on different grid configurations show that 
	the limiter function $H_{3\text{L}}^{(c)}$ maintains the optimal third-order accuracy on smooth profiles 
	and avoids oscillations in case of discontinuous solutions.	
\end{abstract}

\end{frontmatter}
%
%----------------------------------------------------------------------------------------------------------------------------------------------------------------------
\section{Introduction}
%----------------------------------------------------------------------------------------------------------------------------------------------------------------------
%
In the context of finite volume schemes, one of the building blocks for obtaining higher-order accuracy 
is the reconstruction of interface values \cite{LeVeque1992}. There are many different approaches using 
linear and non-linear reconstruction functions. Restricting ourselves to the compact stencil of only the 
cell of interest and its direct neighbors, the best order of accuracy we can obtain is third order. This can 
e.g.~be achieved by constructing a quadratic polynomial whose average over each of the three cells 
of interest needs to agree with the cell average of the solution in these cells.s This reconstruction yields 
third-order accuracy, however, the resulting scheme is linear, causing oscillations at discontinuities 
\cite{Godunov1959}. One possibility to work around this are limiter functions in the MUSCL 
framework \cite{VanLeer1979}. These limiters use three cell mean values per reconstruction 
and are total variation diminishing (TVD), however, they generally yield second-order accuracy, 
see \cite{LeVeque1992} and references therein. One of their major drawbacks is the loss of 
accuracy near extrema \cite{Kriel2017}, also called extrema clipping.

In 1987, Harten et.~ak.~\cite{HartenEngquistOsherChakravarthy1987} presented 
the essentially non-oscillatory (ENO) scheme, further developed by Liu et.~al.~\cite{LiuOsherChan1994} 
to become the weighted essentially non-oscillatory (WENO) scheme. This method enjoys great 
popularity and has been extended by many authors. One of the most-widespread enhancements are 
the smoothness indicators proposed by Jiang and Shu \cite{JiangShu1996} which increase the 
order of accuracy. This scheme will be referred to as WENO-JS. Another development, 
incorporating a global higher-order smoothness indicator was proposed by Borges 
et.~al.~\cite{Borges:2008} and is named WENO-Z. We compare our results to the third-order 
versions of these two methods.

Another approach is the use of non-polynomial reconstructions, such as hyperbolic reconstruction 
schemes, cf.~Marquina \cite{Marquina1994} or local double-logarithmic reconstruction schemes, 
cf.~Artebrant and Schroll \cite{ArtebrantSchroll2006}. Based on their work, {\v C}ada and Torrilhon 
developed a third-order limiter function avoiding the above-mentioned extrema clipping. Our recent 
article \cite{SchmidtmannSeiboldTorrilhon2015} continued this work, introducing a third-order limiter 
function \Hlim. This function contains a decision criterion able to distinguish between smooth 
extrema and discontinuities. This limiter was first developed and tested for one-dimensional 
hyperbolic conservation laws on uniform grids in the context of finite volume methods. Now, the 
aim is to extend the scheme to make it applicable for numerical test cases on non-uniform meshes 
and in two space dimensions.

The paper is structured as follows: in Sec.~\ref{sec:limiter} we recall the formulation of the third-order limiter function $H_{3\text{L}}^{(c)}$ in one space dimension for equidistant grids. Then, the limiter is extended for the use of non-equidistant grids. Sec.~\ref{sec:limiter2D} explains the extension to two space dimensions and the order-fix which allows to maintain high-order accuracy within the flux-splitting framework. Here, we firstly treat Cartesian grids in Sec.~\ref{sec:2DCartesian}, then introduce a parallel adaptive-mesh-refinement (AMR) framework in Sec.~\ref{sec:AMR}, and in Sec.~\ref{sec:2DnonUniform} extend the theory of non-uniform grids to two space dimensions. Numerical results visualizing the theoretical concepts are presented in Sec.~\ref{sec:numerics} and in Sec.~\ref{sec:Conclusion} we draw some conclusions. Finally, more details on the formulation of the limiter function on uniform as well as non-uniform grids can be found in the appendix.
%
%----------------------------------------------------------------------------------------------------------------------------------------------------------------------
\section{Third-Order Limiter in One Space Dimension}\label{sec:limiter}
%----------------------------------------------------------------------------------------------------------------------------------------------------------------------
%
\begin{figure}[t]
	\centering
	\includegraphics[width=0.8\textwidth]{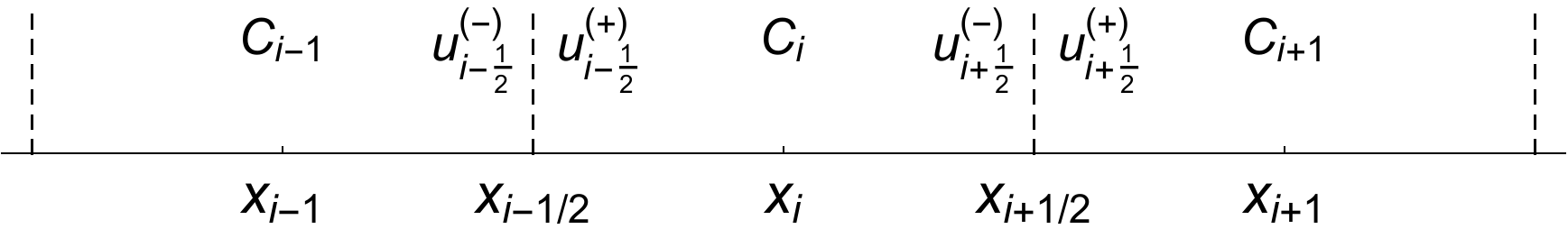}
	\caption{Stencil for reconstruction of interface values on equidistant grids, i.e. $\Delta x_i =x_{i+\frac12} - x_{i-\frac12} \equiv \Delta x\;\forall i$.}
	\label{fig:reconstruction}
\end{figure}
	Achieving high-order accuracy with finite volume schemes requires large stencils. 
	This leads to an increase in communication among grid cells and is undesirable when thinking of 
	parallel codes and boundaries. Therefore, we want to remain on the most compact stencil of three 
	cells in one-dimension and five cells in two space dimensions. This means, the stencil consists of the 
	cell of interest and its 	direct neighbors. 

	For the sake of simplicity, we restrict the theoretical development of the limiter function to 
	one dimensional scalar equations. The transition from the one-dimensional formulation to two-dimensions 
	is obtained via a dimensional splitting. The exact procedure is explained in Sec.~\ref{sec:limiter2D}. 
	Also, the theory easily extends to systems of conservation laws by applying it component-wise.
	
	In the one-dimensional case of Cartesian grids, we divide the domain of interest $\Omega\subset \mathbb{R}$ 
	in non-overlapping cells $C_i=[x_{i-\tfrac12}, x_{i+\tfrac12})$ such that $\Omega=\bigcup_i C_i$. Denote by 
	$x_i$ the cell centers and by $\Delta x_i = x_{i+\frac12} - x_{i-\frac12}$ the size of cell $C_i$. 
	Fig.~\ref{fig:reconstruction} depicts the here-introduced notation for the equidistant case 
	$\Delta x_i \equiv \Delta x\;\forall i$. 
	\\	
	\\	
	In this work we are interested in hyperbolic conservation law of the form
	\begin{align}
		\partial_t u(\vec{x},t) + \nabla\cdot f(u(\vec{x},t)) = 0
		\label{eq:consLaw}
	\end{align}
	with suitable initial conditions $u(\vec{x}, 0)=u_0(\vec{x})$, $\vec{x}=x$ in 1D and $\vec{x}=(x,y)^T$ in 2D, 
	respectively. To avoid boundary effects, we impose periodic boundary conditions. Integrating 
	Eq.~\eqref{eq:consLaw} in one space dimension over cell $C_i$ and dividing by the cell width $\Delta x_i$ 
	yields an exact update formula for the cell mean values $\tfrac{1}{\Delta x_i}\int_{C_i} u(x,t) dx$. This 
	formulation however requires the exact solution of Riemann problems at each cell boundary 
	\cite{LeVeque1992}. To avoid this costly procedure, approximate Riemann solvers are incorporated and the 
	exact cell mean values are approximated by $\bar{u}_i$. The update formula for $\bar{u}_i$ is then given by 
	the so-called semi-discrete scheme
	\begin{align}
		\frac{d\,\bar{u}_i}{d\,t} = -\frac{1}{\Delta x_i} \left( \hat{f}_{i+\frac12} - \hat{f}_{i-\frac12} \right),
		\label{eq:FVsemidiscrete}
	\end{align}
	with numerical flux functions $\hat{f}_{i+\frac12} = \hat{f}(u^{(-)}_{i+\frac12}, u^{(+)}_{i+\frac12})$ and 
	$\hat{f}_{i-\frac12} = \hat{f}(u^{(-)}_{i-\frac12}, u^{(+)}_{i-\frac12})$. These functions take as input the left and 
	right limiting values at the cell interface, see Fig.~\ref{fig:reconstruction}. One could simply insert the left and 
	right cell mean values, however, this only yields a first-order accurate scheme \cite{LeVeque1992}. In order to 
	achieve higher-order accuracy, one way is to use reconstructions for the interface values. As described in more 
	detail in \cite{SchmidtmannSeiboldTorrilhon2015}, the reconstructed interface values of cell $i$ can be written 
	in the general form
	\begin{subequations}
		\begin{align}
	 	  	u^{(-)}_{i+\frac12} &= \bar u_i + \tfrac{1}{2}\; H(\delta_{i-\frac{1}{2}},\delta_{i+\frac{1}{2}}) = M(\bar u_{i-1}, \bar u_{i}, \bar u_{i+1}), \\
  		  	u^{(+)}_{i-\frac12} &= \bar u_i - \tfrac{1}{2}\; H(\delta_{i+\frac{1}{2}},\delta_{i-\frac{1}{2}}) = P(\bar u_{i-1}, \bar u_{i}, \bar u_{i+1}).
		\end{align}
	    	\label{eq:reconstruction}
  	\end{subequations}
	\hspace{-0.3cm} The function $H$ fully determines the way limiting is performed and thus the order of 
	accuracy of the resulting scheme. The undivided differences between neighboring cells are denoted by
	\begin{align}
		\begin{split}
			\delta_{i-\frac{1}{2}} = \bar u_i-\bar u_{i-1}\\
			\delta_{i+\frac{1}{2}} = \bar u_{i+1}-\bar u_i.
		\end{split}
		\label{eq:deltas}
	\end{align}
	\paragraph{Remark:}
	The standard form for reconstructions, e.g.~found in \cite{LeVeque1992} reads
	\begin{subequations}
		\begin{align}
	 	  	u^{(-)}_{i+\frac12} &= \bar u_i + \tfrac{1}{2}\; \phi(\theta_i)\delta_{i-\frac{1}{2}}, \\
  		  	u^{(+)}_{i-\frac12} &= \bar u_i - \tfrac{1}{2}\; \phi(\theta_i^{-1})\delta_{i+\frac{1}{2}}
		\end{align}
	    	\label{eq:reconstructionStandard}
  	\end{subequations}
	\hspace{-0.3cm} with the ratio of consecutive gradients 
	$\theta_i = \delta_{i-\frac{1}{2}}/\delta_{i+\frac{1}{2}}$ which acts as 
	a smoothness indicator and a monovariant limiter function $\phi$. This form, introduced 
	by \cite{Dubois1990, CadaTorrilhon2009} relates to Eq.~\eqref{eq:reconstruction} via 
	$\phi(\theta_i)\delta_{i-\frac{1}{2}} = H(\delta_{i-\frac{1}{2}},\delta_{i+\frac{1}{2}})$.
	More details on the two-variate form and its advantages can be found in 
	\cite{SchmidtmannSeiboldTorrilhon2015}.
%
%----------------------------------------------------------------------------------------------------------------------------------------------------------------------
\subsection{Formulation for Equidistant Grids}\label{sec:equidistant}
%----------------------------------------------------------------------------------------------------------------------------------------------------------------------
%
	In this section we will shortly recall the formulation of third-order limiter functions for equidistant grids 
	developed in \cite{SchmidtmannSeiboldTorrilhon2015}.\\
	Starting with a quadratic ansatz function, evaluated on $(x_{i-1}, \bar u_{i-1}), (x_i, \bar u_i), 
	(x_{i+1}, \bar u_{i+1})$, we can obtain an unlimited reconstruction formulation which yields a 
	third-order accurate scheme. Rewriting the polynomial reconstruction in the form 
	\eqref{eq:reconstruction} yields
	\begin{align}
		&u^{(\mp)}_{i\pm\frac12} = \bar u_i \pm \frac{1}{2}\; \frac{\delta_{i+\frac{1}{2}} + \delta_{i\pm\frac{1}{2}} + \delta_{i-\frac{1}{2}}}{3}\\
		\intertext{leading to the function}
		&H_{3}(\delta_{i-\frac{1}{2}}, \delta_{i+\frac{1}{2}}) :=\frac{1}{3}\left(2 \delta_{i+\frac{1}{2}} + \delta_{i-\frac{1}{2}}\right).
		\label{eq:3rdOrder}			
	\end{align}
	For purely smooth functions, the full (unlimited) third-order reconstruction shows good results. However, 
	for solutions containing discontinuities, spurious oscillations develop since a linear higher-order method 
	is not monotonicity preserving \cite{LeVeque1992}. 
	Therefore, we need to apply non-linear reconstruction functions. In \cite{SchmidtmannSeiboldTorrilhon2015} 
	we recently constructed a limiter function, called $H_{3\text{L}}$, based on a double logarithmic ansatz 
	function, first proposed by Artebrant and Schroll \cite{ArtebrantSchroll2006} and further developed by 
	{\v{C}}ada and Torillhon \cite{CadaTorrilhon2009}. 
	Furthermore, we designed the combined limiter function $H_{3\text{L}}^{(c)}$ which 
	includes a decision criterion $\eta = \eta(\delta_{i-\frac{1}{2}},\delta_{i+\frac{1}{2}})$, 
	able to distinguish between smooth extrema and discontinuities. 
	The combined limiter applies the full third-order reconstruction $H_3$ on parts which are classified as smooth 
	and switches to the limited function $H_{3\text{L}}$ if $\eta$ indicates large gradients. It reads
	\begin{align}     
		H_{3\text{L}}^{(c)} (\delta_{i-\frac{1}{2}},\delta_{i+\frac{1}{2}}) := 
			\begin{cases}
				H_{3}(\delta_{i-\frac{1}{2}},\delta_{i+\frac{1}{2}}) \quad\qquad &\text{if}\;\eta(\delta_{i-\frac{1}{2}},\delta_{i+\frac{1}{2}}) <1  \\
		  	  	H_{3\text{L}}(\delta_{i-\frac{1}{2}},\delta_{i+\frac{1}{2}})\quad\; &\text{if}\;\eta(\delta_{i-\frac{1}{2}},\delta_{i+\frac{1}{2}}) \geq 1.
		  	\end{cases}
		  \label{eq:H3Lc}	
	\end{align}	
	The formulation for $H_{3\text{L}}$ and the decision criterion $\eta$, as well as more details on 
	$H_{3\text{L}}^{(c)}$ are given in \ref{sec:AppendixA} and in \cite{SchmidtmannSeiboldTorrilhon2015}.
	
%	The reconstrution obtained by Eq.~\eqref{eq:reconstruction}, using the full-third-order 
%	reconstruction $H_3(\cdot,\cdot)$, and the limiter function $H_{3\text{L}}^{(c)}$, respectively, 
%	is homogeneous, translational invariant, and symmetric, see \cite{SchmidtmannSeiboldTorrilhon2015} 
%	for more details.
%
%----------------------------------------------------------------------------------------------------------------------------------------------------------------------
\subsection{Formulation for Non-Equidistant Grids}\label{sec:nonUniform}
%----------------------------------------------------------------------------------------------------------------------------------------------------------------------
%
\begin{figure}[t]
	\centering
	\includegraphics[width=0.8\textwidth]{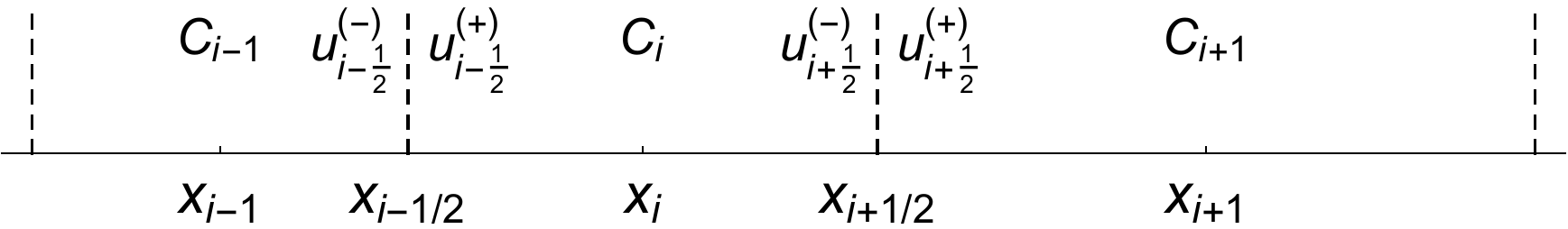}
	\caption{Stencil for reconstruction of interface values in the non-equidistant case.}
	\label{fig:reconstructionNeq}
\end{figure}
	For general grids, the size of cell $C_i$, denoted by $\Delta x_i$ is not uniform for 
	all cells, i.e. $\Delta x_i \neq \Delta x\;\forall i$, see Fig.~\ref{fig:reconstructionNeq}. 
	In this case, the definition of the undivided differences $\delta _{i\pm\frac{1}{2}}$ 
	Eq.~\eqref{eq:deltas} is not meaningful anymore and new concepts need to be developed. 
	
	Starting again with the full third-order reconstruction, consider a quadratic polynomial 
	$p_i(x)$ in cell $i$ that has to maintain the cell averages in the three cells 
	$C_{i+\ell},\; \ell\in\{-1, 0, 1\}$. This polynomial is then evaluated at the cell boundaries 
	$x_{i\pm \frac12}$ and yields the reconstructed cell interface values
	\begin{subequations}
		\begin{align}			
			u_{i+\frac12}^{(-)} &= p_i(x_{i+\frac{1}{2}}) \stackrel{!}{=} \bar u_i + \frac{1}{2} H_{3,\text{neq}}\\
		u_{i-\frac12}^{(+)} &= p_i(x_{i-\frac{1}{2}}) \stackrel{!}{=} \bar u_i - \frac{1}{2}H_{3,\text{neq}}.
		\end{align}
		\label{eq:reconstructionNeq}		
	\end{subequations}
	Even though this procedure is similar to the full third-order reconstruction on uniform grids, 
	Eq.~\eqref{eq:3rdOrder}, the reconstruction function $H_{3,\text{neq}}$ differs from $H_{3}$ since the 
	different cell sizes need to be taken into account. The full (unlimited) third-order reconstruction function 
	reads
	\begin{subequations}
		\begin{align}
			&H_{3,\text{neq}}\left(\delta _{i-\frac{1}{2}},\delta _{i+\frac{1}{2}},	\Delta x_i, \Delta x_{i-1}, \Delta x_{i+1}\right) = \frac{\Delta x_i}{\Delta_i} \frac{1}{3}\left(2 \frac{\Delta_{i-\frac12}}{\Delta_{i+\frac12}}\delta _{i+\frac{1}{2}} + \frac{\Delta x_{i+1}}{\Delta_{i-\frac12}}\delta _{i-\frac{1}{2}}\right)
			\label{eq:3rdOrderNeqH}
			\intertext{with the abbreviations}
			&\Delta_i = \frac{\Delta x_{i-1}+\Delta x_i+\Delta x_{i+1}}{3}, \;\,\Delta_{i-\frac12} = 
			\frac{\Delta x_{i-1}+\Delta x_{i}}{2},\;\,\Delta_{i+\frac12} = \frac{\Delta x_i+\Delta x_{i+1}}{2}.		
			\label{eq:3rdOrderNeqabbreviations}			
		\end{align}
		\label{eq:3rdOrderNeq}
	\end{subequations}
	As in the equidistant case, the reconstructed interface values, Eq.~\eqref{eq:reconstructionNeq}, can 
	be compactly refomulated as
	\begin{subequations}
		\begin{align}		
			u^{(\mp)}_{i\pm\frac12} = \bar u_i \pm 
			\frac{\Delta x_i}{2}\; 
			\frac{\Delta x_{i-1}\widetilde \delta _{i+\frac{1}{2}} + \Delta x_i \widetilde\delta_{i\pm\frac{1}{2}} 
			+\Delta x_{i+1}\widetilde	 \delta_{i-\frac{1}{2}}}{\Delta x_{i-1}+\Delta x_i+\Delta x_{i+1}}		
			\label{eq:uInterfaceRec3rdOrderNeq}
			\intertext{with}
			\widetilde \delta _{i-\frac{1}{2}} = \frac{\delta _{i-\frac{1}{2}}}{\Delta_{i-\frac12}},\qquad
			\widetilde \delta _{i+\frac{1}{2}} = \frac{\delta _{i+\frac{1}{2}}}{\Delta_{i+\frac12}}.
			\label{eq:deltasNeq}		
		\end{align}
		\label{eq:3rdOneqCompact}
	\end{subequations}

	It can easily be seen that for equidistant grids, i.e.~$\Delta x_{i-1}=\Delta x_{i}=\Delta x_{i+1}\equiv \Delta x$,
	the abbreviated terms reduce to $\Delta_i = \Delta_{i-\frac12} = \Delta_{i+\frac12} = \Delta x$ and 
	therefore, the formulas for $H_{3,\text{neq}}$ and $H_{3}$ match, as expected.

	Eq.~\eqref{eq:3rdOrderNeq} and \eqref{eq:3rdOneqCompact} indicate that for non-equidistant meshes, 
	the equivalent of the undivided differences $\delta_{i\pm\frac{1}{2}}$ are the scaled slopes
	\begin{align}
		u_{i+\frac12}^{(-)} : \qquad
		\begin{split}
		\delta _{i-\frac{1}{2}} \; &\rightarrow \; \Delta x_{i+1} \frac{\delta _{i-\frac{1}{2}}}{\Delta_{i-\frac12}} 
		= \Delta x_{i+1} \widetilde \delta _{i-\frac{1}{2}} \\
		\delta _{i+\frac{1}{2}} \; &\rightarrow \; \Delta_{i-\frac12}\;\frac{\delta _{i+\frac{1}{2}} }{\Delta_{i+\frac12}}
		= \frac{\Delta x_i}{2}\;\widetilde \delta _{i+\frac{1}{2}} + \frac{\Delta x_{i-1}}{2}\;\widetilde \delta _{i+\frac{1}{2}}
		\end{split}
		\label{eq:deltasNeqRight}
	\end{align}
	for the reconstruction of the right interface of cell $C_i$ and
	\begin{align}
		u_{i-\frac12}^{(+)}: \qquad 
		\begin{split}
		\delta _{i-\frac{1}{2}} \; &\rightarrow \; \Delta_{i+\frac12}\;\frac{\delta_{i-\frac{1}{2}}}{\Delta_{i-\frac12}}
		=  \frac{\Delta x_i}{2}\;\widetilde \delta _{i-\frac{1}{2}} + \frac{\Delta x_{i+1}}{2}\;\widetilde \delta _{i-\frac{1}{2}} \\
		\delta _{i+\frac{1}{2}} \; &\rightarrow \; \Delta x_{i-1}\frac{\delta_{i+\frac{1}{2}}}{\Delta_{i+\frac12}}	
		=\Delta x_{i-1} \widetilde \delta _{i+\frac{1}{2}}
		\end{split}
		\label{eq:deltasNeqLeft}		
	\end{align}
	for the reconstruction of the left interface of cell $C_i$. 
	These expressions resemble the smoothness indicators introduced by 
	Jiang and Shu \cite{JiangShu1996}, which are given by 
	$\Delta x_i \widetilde \delta _{i-\frac{1}{2}}$ and 
	$\Delta x_i \widetilde \delta _{i+\frac{1}{2}}$.
%	Let us take a closer look at the developed smoothness indicators \eqref{eq:deltasNeqRight} and \eqref{eq:deltasNeqLeft}. For the reconstruction of $u^{(-)}_{i+\tfrac{1}{2}}$, 
	
	In order to generalize the third-order limiter function developed in \cite{SchmidtmannSeiboldTorrilhon2015}, 
	we replace the undivided differences as mentioned above to obtain the reconstructions
	\begin{subequations}
		\begin{align}			
			u_{i+\frac12}^{(-)} &= \bar u_i + \frac{1}{2} H^{(c)}_{3L}\left( \Delta x_{i+1} \widetilde \delta _{i-\frac{1}{2}}, \Delta_{i-\frac12}\;\widetilde \delta _{i+\frac{1}{2}} \right)\\
		u_{i-\frac12}^{(+)} &= \bar u_i - \frac{1}{2}H^{(c)}_{3L}\left( \Delta x_{i-1}\widetilde \delta _{i+\frac{1}{2}}, \Delta_{i+\frac12}\;\widetilde \delta _{i-\frac{1}{2}} \right).
		\end{align}
		\label{eq:limiterReconstructionNeq}
	\end{subequations}
	with the limiter function $H^{(c)}_{3L}$ \eqref{eq:H3Lc} described in Sec.~\ref{sec:equidistant}. 
	The non-equidistant version of the limiter function can be defined as a function $H^{(c)}_{3L,\text{neq}}$, 
	given by
	\begin{align}
	\begin{split}
		H^{(c)}_{3L,\text{neq}}\left(\delta _{i-\frac{1}{2}},\delta _{i+\frac{1}{2}},	\Delta x_i, \Delta x_{i-1}, \Delta x_{i+1}\right) &= H^{(c)}_{3L}\left( \Delta x_{i+1} \frac{\delta _{i-\frac{1}{2}}}{\Delta_{i-\frac12}}, \Delta_{i-\frac12}\;\frac{\delta _{i+\frac{1}{2}} }{\Delta_{i+\frac12}} \right) \\
		H^{(c)}_{3L,\text{neq}}\left(\delta _{i+\frac{1}{2}},\delta _{i-\frac{1}{2}},\Delta x_i, \Delta x_{i+1}, \Delta x_{i-1}\right) &= H^{(c)}_{3L}\left( \Delta x_{i-1}\frac{\delta_{i+\frac{1}{2}}}{\Delta_{i+\frac12}}, \Delta_{i+\frac12}\;\frac{\delta_{i-\frac{1}{2}}}{\Delta_{i-\frac12}} \right).
	\end{split}	
	\end{align}
	The decision criterion $\eta$ \eqref{eq:eta} for non-uniform meshes reads
	\begin{align}
		\eta(\delta_1,\delta_2) = \frac{\sqrt{\delta_1^2+\delta_2^2}}{\sqrt{\frac{5}{2}}\,\alpha\,dx^2}, 
	\end{align}
	where $\delta_1,\delta_2$ are the same input arguments as for $H^{(c)}_{3L,\text{neq}}$, see 
	Eq.~\eqref{eq:H3Lc} and $dx$ is the average mesh size, $dx = (\sum_i \Delta x_i)/\# cells$.
%
%----------------------------------------------------------------------------------------------------------------------------------------------------------------------
\section{Third-Order Limiter in Two Space Dimensions}\label{sec:limiter2D}
%----------------------------------------------------------------------------------------------------------------------------------------------------------------------
%
	In this section we extend the third-order limiter function to two space dimensions covering 
	three core areas. First, we discuss how to apply the scheme on uniform Cartesian grids. 
	Then we extend the method to adaptively refined grids. 
	The last part of this section explains how the method can be used on rectangular grids 
	which are non-uniform in $x$- and $y$-direction. 
%
%----------------------------------------------------------------------------------------------------------------------------------------------------------------------
\subsection{Formulation for 2D Cartesian Grids}\label{sec:2DCartesian}
%----------------------------------------------------------------------------------------------------------------------------------------------------------------------
%
	In two space dimensions, the domain of interest $\Omega$ is divided into non-overlapping 
	cells 	$C_{i,j}=[x_{i-\frac12}, x_{i+\frac12})\times [y_{j-\frac12}, y_{j+\frac12})$ such that
	$\Omega=\bigcup_{i,j} C_{i,j}$. 
	Denote by $(x_{i},y_{j})$ the cell center of cell $C_{i,j}$.
	The mesh width is given by $\Delta x_{i} = x_{i+\frac12} - x_{i-\frac12}$ and 
	$\Delta y_{j} = y_{j+\frac12} - y_{j-\frac12}$.
	Furthermore we denote by $\tilde{\bar{u}}_{i,j}$ the cell--averaged value over cell $C_{i,j}$
	and by $\tilde{u}_{i+\frac12,j}$ and $\bar{u}_{i,j+\frac12}$ the interface--averaged
	values over the corresponding interface. The tilde notation $\tilde{\cdot}$ denotes the average
	in $y$-direction and bar $\bar{\cdot}$ denotes the average in $x$-direction as in the 
	one-dimensional case.
	
	Integrating a hyperbolic conservation law of the form
	\begin{align}
	  \partial_t u(x,y,t) + \partial_x f(u(x,y,t)) + \partial_y g(u(x,y,t)) = 0 
	  \label{eq:consLaw2d}
	\end{align}
	over cell $C_{i,j}$ and dividing by the cell area $\Delta x_i\Delta y_j$ yields the two-dimensional 
	semi-discrete flux-differencing finite volume scheme \cite{LeVeque2002}
	\begin{align}
	  \frac{d\,\tilde{\bar{u}}_{i,j}}{d\,t} = 
	  -\frac{1}{\Delta x_{i}} \left(\tilde {\hat{f}}_{i+\frac12,j} - \tilde{\hat{f}}_{i-\frac12,j} \right) 
	  -\frac{1}{\Delta y_{j}} \left( \bar {\hat{g}}_{i,j+\frac12} - \bar {\hat{g}}_{i,j-\frac12} \right).
	\end{align}
	Here, the numerical flux functions ${\hat{f}}_{i+\frac12,j}$ and  ${\hat{g}}_{i,j+\frac12}$ are
	approximations to averages of the flux across the corresponding interface \cite{LeVeque2002}
	\begin{subequations}
	  \begin{align}
	   \tilde {\hat{f}}_{i+\frac12,j}&\approx 
	                       \frac{1}{\Delta y_{j}}\int_{y_{j-\frac{1}{2}}}^{y_{j+\frac{1}{2}}} f(u(x_{i+\frac{1}{2}},y,t)) dy,\\
	   \bar {\hat{g}}_{i,j+\frac12} &\approx
	                        \frac{1}{\Delta x_{i}}  \int_{x_{i-\frac{1}{2}}}^{x_{i+\frac{1}{2}}} g(u(x,y_{j+\frac{1}{2}},t)) dx.  
	  \end{align}
	\end{subequations}
	Similar to other finite volumen methods, applying the scheme described in Sec.~\ref{sec:limiter} in 
	a dimension--by--dimension fashion results in a second order scheme,
	see e.g.~\cite{Shu:2009}\cite{Zhang:2011}\cite{BuchmuellerHelzel2014}.
	In order to remain third-order accurate, we apply the fourth order transformation proposed 
	by Buchm\"uller and Helzel 	\cite{BuchmuellerHelzel2014}. Incorporating this so-called order-fix, 
	the scheme can be summarized as follows.
	\begin{enumerate}
	\item Compute the averaged values of the conserved quantities at the cell interfaces in 
		the interior of cell $C_{i,j}$ for all $i,j$ using the one-dimensional limiter functions described 
		in Sec.~\ref{sec:equidistant}
		 \begin{subequations}
			\label{eq:reconstruction2d}
		    \begin{align}
		    \begin{split}
%			    \tilde  u^{(-)}_{i+\frac12,j} &= \tilde{\bar u}_{i,j} + \tfrac{1}{2}\; H(\delta_{i-\frac{1}{2},j},\delta_{i+\frac{1}{2},j}), 
%			     &\tilde u^{(+)}_{i-\frac12,j} &= \tilde{\bar u}_{i,j} - \tfrac{1}{2}\; H(\delta_{i+\frac{1}{2},j},\delta_{i-\frac{1}{2},j}), 	  
%			    \bar  u^{(-)}_{i,j+\frac12} &= \tilde{\bar u}_{i,j} + \tfrac{1}{2}\; H(\delta_{i,j-\frac{1}{2}},\delta_{i,j+\frac{1}{2}}), 
%			     &\bar u^{(+)}_{i,j-\frac12} &= \tilde{\bar u}_{i,j} - \tfrac{1}{2}\; H(\delta_{i,j+\frac{1}{2}},\delta_{i,j-\frac{1}{2}}).
    \tilde  u^{(-)}_{i+\frac12,j} &= \tilde{\bar u}_{i,j} + \tfrac{1}{2}\; H(\delta_{i-\frac{1}{2},j},\delta_{i+\frac{1}{2},j}), \\
		     \tilde u^{(+)}_{i-\frac12,j} &= \tilde{\bar u}_{i,j} - \tfrac{1}{2}\; H(\delta_{i+\frac{1}{2},j},\delta_{i-\frac{1}{2},j}), \\
		    \bar  u^{(-)}_{i,j+\frac12} &= \tilde{\bar u}_{i,j} + \tfrac{1}{2}\; H(\delta_{i,j-\frac{1}{2}},\delta_{i,j+\frac{1}{2}}), \\
		     \bar u^{(+)}_{i,j-\frac12} &= \tilde{\bar u}_{i,j} - \tfrac{1}{2}\; H(\delta_{i,j+\frac{1}{2}},\delta_{i,j-\frac{1}{2}}).			     
		    \end{split}
		    \end{align}
	\end{subequations}
  	Here, the reconstruction function $H$ can be the unlimited third-order reconstruction $H_3$, the limiter 
  	function $H_{3L}^{(c)}$ or any other third-order limiter fitting the setting.
  		The undivided differences in two-dimensions, $\delta_{i\pm\frac{1}{2},j}, \delta_{i, j\pm\frac{1}{2}}$, 
  		are defined similarly to their one-dimensional equivalents, Eq.~\eqref{eq:deltas},
  		\begin{align}
  			\begin{split}
			\delta_{i-\frac{1}{2},j} &= \tilde{\bar u}_{i,j} -\tilde{\bar u}_{i-1,j}\\
			\delta_{i+\frac{1}{2},j} &= \tilde{\bar u}_{i+1, j} -\tilde{\bar u}_{i,j}\\
			\delta_{i, j-\frac{1}{2}} &= \tilde{\bar u}_{i,j} -\tilde{\bar u}_{i,j-1}\\
			\delta_{i, j+\frac{1}{2}} &= \tilde{\bar u}_{i,j+1} -\tilde{\bar u}_{i,j}.
		\end{split}
		\label{eq:deltas2D}
  		\end{align}
	\item Compute point values of the conserved quantities at the center of each cell
	  interface, i.e.\ compute 
	  \begin{equation}\label{eqn:ave2point_4th}
	    \begin{split}
	      u_{i+\frac{1}{2},j}^{(\pm)} & = \tilde u_{i+\frac{1}{2},j}^{(\pm)} - \frac{1}{24}
	      \left( \tilde u_{i+\frac{1}{2},j-1}^{(\pm)} - 2 \tilde u_{i+\frac{1}{2},j}^{(\pm)} +
	        \tilde u_{i+\frac{1}{2},j+1}^{(\pm)} \right), \\
	      u_{i,j+\frac{1}{2}}^{(\pm)} & = \bar u_{i,j+\frac{1}{2}}^{(\pm)} - \frac{1}{24}
	      \left( \bar u_{i-1,j+\frac{1}{2}}^{(\pm)} - 2 \bar u_{i,j+\frac{1}{2}}^{(\pm)} +
	        \bar u_{i+1,j+\frac{1}{2}}^{(\pm)} \right).
	    \end{split}
	  \end{equation}
	\item Compute fluxes at the center of the cell interfaces using the computed point values 
		and a consistent numerical flux function, i.e.\ 
		  \begin{align}
			  {\hat{f}}_{i+\frac12,j}={\hat{f}}(u^{(-)}_{i+\frac12,j}, u^{(+)}_{i+\frac12,j}), \quad
			  {\hat{g}}_{i,j+\frac12}={\hat{g}}(u^{(-)}_{i,j+\frac12}, u^{(+)}_{i,j+\frac12}).    
		  \end{align}
	\item Compute averaged values of the numerical flux function, i.e.\ compute 
		\begin{equation} \label{eqn:flux-4th}
		    \begin{split}
		       \tilde{\hat{f}}_{i+\frac12,j} & = \hat{f}_{i+\frac{1}{2},j} + \frac{1}{24} \left(
		        \hat{f}_{i+\frac{1}{2},j-1} - 2 \hat{f}_{i+\frac{1}{2},j} +
		        \hat{f}_{i+\frac{1}{2},j+1} \right), \\
		     \bar {\hat{g}} _{i,j+\frac{1}{2}} & = \hat{g}_{i,j+\frac{1}{2}} + \frac{1}{24} \left(
		        \hat{g}_{i-1,j+\frac{1}{2}} - 2 \hat{g}_{i,j+\frac{1}{2}} +
		        \hat{g}_{i+1,j+\frac{1}{2}} \right).
		    \end{split}
	  \end{equation}
	\item Use a high--order accurate Runge--Kutta method for the update in time. 
	In this work, we use the strong stability preserving third--order Runge--Kutta method 
	described by Gottlieb et.~al.~\cite{GottliebShuTadmor2001}.
	\end{enumerate}

This procedure is quite robust even when discontinuities are present. 
Nevertheless in some situations an unphysical state may be created, therefor we apply a
simple limiting as suggested by Buchm\"uller et.~al.~\cite{Buchmueller:2016}. 
Details on this limiting procedure can be found in the original paper.

Step 1. comprises the reconstruction function $H(\cdot,\cdot)$ that has been described in 
Sec.~\ref{sec:limiter}. For purely smooth solutions, the full third-order reconstruction 
$H_3$ can be used in this step. However, when discontinuities are present, the limiter function 
$H_{3\text{L}}^{(c)}$, Eq.~\eqref{eq:H3Lc}, is more advisable since oscillations are prevented. 
In principle, in the dimension-splitting approach described above, $H_{3\text{L}}^{(c)}$ can 
be applied in the same manner as in one dimension. The decision criterion 
$\eta$ (see Appendix \ref{sec:AppendixA}, Eq.~\eqref{eq:eta}) can also be used without 
any changes in the splitting approach. Only the definition of the radius of the asymptotic 
region, $\alpha$ (see Appendix \ref{sec:AppendixA}, Eq.~\eqref{eq:alpha}) has to be adapted. 
It is defined as 
\begin{align}
	\alpha = \max\limits_{(x,y)\in\Omega\backslash \Omega_d}|\Delta u_0(x,y)|
	\label{eq:alpha2D}
\end{align}
in the two-dimensional scheme. Again, $\Omega$ is the domain of interest and $\Omega_d\subset \Omega$ the subset containing discontinuities.
%
%----------------------------------------------------------------------------------------------------------------------------------------------------------------------
\subsection{Adaptive Mesh Refinement (AMR)}\label{sec:AMR}
\begin{figure}[t]
  \centering
  \subfloat[][Hierarchy of grid blocks.]{% 
    \begin{tikzpicture}[scale=.99*0.72]
      \draw[fill=black!10!white] (0,0) grid [xstep=2.,ystep=2.] (8,8) rectangle (0,0);
      \draw[fill=black!20!white] (0,0) grid [xstep=1.,ystep=1.] (6,6)  rectangle (0,0);
      \draw[fill=black!30!white] (2,2) grid [xstep=.5,ystep=.5] (4,5)  rectangle (2,2);
      \draw[fill=black!30!white] (4,3) grid [xstep=.5,ystep=.5] (5,5)  rectangle (4,3);
      \draw[fill=black!40!white] (2.5,3.) grid [xstep=.25,ystep=.25] (3.5,4.5)  rectangle (2.5,3.);
      \draw[fill=black!40!white] (3.5,3.5) grid [xstep=.25,ystep=.25] (4,4.5)  rectangle (3.5,3.5);
    \end{tikzpicture}
    \label{fig:racoon-grid-struc}
  }\;
  \subfloat[][Zoom of a single block.]{%  
    \begin{tikzpicture}[scale=0.79*0.72]
      \draw[draw=black!50!white, fill=black!10!white]
      (0,0) grid [xstep=0.5,ystep=0.5] (10,10) rectangle (0,0);
      \draw[ fill=white] (1,1) grid [xstep=0.5,ystep=0.5] (9,9) rectangle (1,1);
    \end{tikzpicture}
    \label{fig:racoon-grid-block}
  }
  \caption{\protect\subref{fig:racoon-grid-struc} Hierarchy of grid
    blocks of level $2-5$. \protect \subref{fig:racoon-grid-block} A single block
    consisting of $16\times 16$ grid cells, and a layer of two ghost cells.}
\end{figure}
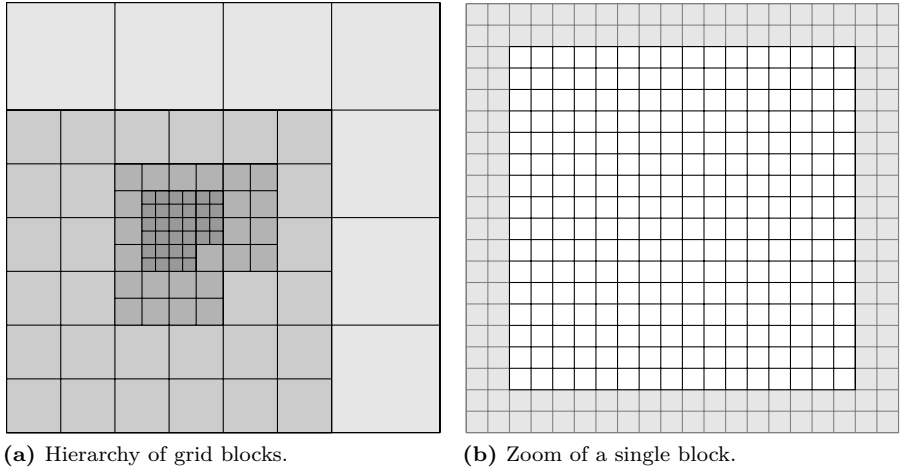

%----------------------------------------------------------------------------------------------------------------------------------------------------------------------
% 
For computations in two dimensions, we use the parallel AMR framework
{\em  Racoon} developed by Dreher and Grauer \cite{Dreher:2005}. 
Both, the grid adaptivity and the parallelization are based on a block--structure. 
In a $2$-dimensional space, a grid of level $\ell$ consists of $(2^2)^\ell = 4^\ell$ blocks.
Computations can be performed simultaneously  on each block and due to the Cartesian
grid structure within each block, we can simply apply the method described in 
Sec.~\ref{sec:2DCartesian}. 
A typical block is illustrated in Fig.~\ref{fig:racoon-grid-block}. 
The cells in the gray region are ghost cells needed for the communication between the
blocks. 
For refinement a block of level $\ell$ is replaced by $2^d$ blocks of level $\ell+1$.
These blocks may then be further refined until the maximum refinement level is reached. 
There are three reasons for refinement.

\FloatBarrier

\begin{itemize}
\item Some refinement criteria is met. 
  Here we compute 
  \begin{equation}\label{eqn:refinement}
    \delta = \frac{|q_{i-1,j}  - 2 q_{i,j} + q_{i+1,j}|
    +|q_{i,j-1} - 2 q_{i,j} + q_{i,j+1}| }{|q_{i,j}| \Delta x  \Delta y}
  \end{equation}
  for each cell. 
  If $\delta$ is bigger than a predefined threshold $\delta_0$, the cell is marked for
  refinement and therefore the block will be refined.
\item Neighbouring blocks are may also refined, so that after refinement the region with
  the marked cell is surrounded by fine blocks. 
  In 2D for example, if a cell in the upper left part of a block is marked for refinement,
  then the upper block, the block on the left-hand side, and the block in the upper left
  diagonal direction will be refined as well.
\item Finally the grid needs to be properly nested, that is the level of neighboring
  blocks is not allowed to differ by more than one. 
  Which may lead to further refinement.
\end{itemize}
In Fig.~\ref{fig:racoon-grid-struc}, a typical block structure is illustrated. The blocks
in this figure are of level $2-5$, where light gray corresponds to 2 and increasing 
darkness corresponds to increasing refinement level.

As mentioned above, ghost cells are used for communication between blocks and need to be
updated in every stage of the time stepping scheme.
In most cases this means simply copying the cell--averaged data from the neighboring
block. 
To transfer data from a fine block to a coarse block, the values of the corresponding
cells are averaged. 
Values for the fine block are created by polynomial reconstruction using data of the
coarse block. 
The same is procedure is applied when a block is refined or coarsen again, see 
\cite{BuchmuellerHelzel2014} for more details.  
%
%----------------------------------------------------------------------------------------------------------------------------------------------------------------------
\subsection{Non-Uniform Rectangular 2D Grids}\label{sec:2DnonUniform}
%----------------------------------------------------------------------------------------------------------------------------------------------------------------------
%
In this section we consider non-uniform two-dimensional meshes.
The Cartesian grid cells are transformed into non-uniform cells in $x$- and $y$-direction by adding 
a perturbation to the cell centers $(x_i, y_j)$. In this work we used the transformation 
\begin{align}
	\begin{split}
	&x_i\quad\rightarrow\quad x_i + \delta x \sin(c_x \pi x_i), \\
	&y_j\quad\rightarrow\quad  y_j + \delta y \sin (c_y\pi  y_j)
	\end{split}
	\label{eq:transformationToNonUniform}
\end{align}
with the constants $\delta x, c_x, \delta_y, c_y$, which determine the structure of the mesh.
This procedure yields rectangles that are still aligned with the $x$- and $y$-axes but exhibit different 
cell sizes.
\\
\\
In order to apply the numerical schemes presented in the first sections of this paper, we need to adapt 
the schemes as follows: The general structure of the numerical algorithm for two-dimensional Cartesian 
grids, introduced in Sec.~\ref{sec:2DCartesian}, remains the same. We only need to adapt step 1, i.e. 
the reconstruction of the interface values. The first step of the algorithm for non-uniform grids reads
 \begin{enumerate}
	\item Compute the averaged values of the conserved quantities at the cell interfaces in 
		the interior of cell $C_{i,j}$ for all $i,j$ using the one-dimensional limiter functions for non-uniform 
		grids, described 	in Sec.~\ref{sec:nonUniform}.
		 \begin{subequations}
			\label{eq:reconstruction2dneq}
		    \begin{align}
		    \begin{split}
		\tilde  u^{(-)}_{i+\frac12,j} &= \tilde{\bar u}_{i,j} + \tfrac{1}{2}\; H(\delta_{i-\frac{1}{2},j},\delta_{i+\frac{1}{2},j},\Delta x_i, \Delta x_{i-1}, \Delta x_{i+1}), \\
		     \tilde u^{(+)}_{i-\frac12,j} &= \tilde{\bar u}_{i,j} - \tfrac{1}{2}\; H(\delta_{i+\frac{1}{2},j},\delta_{i-\frac{1}{2},j}, \Delta x_i, \Delta x_{i+1}, \Delta x_{i-1}), \\
		    \bar  u^{(-)}_{i,j+\frac12} &= \tilde{\bar u}_{i,j} + \tfrac{1}{2}\; H(\delta_{i,j-\frac{1}{2}},\delta_{i,j+\frac{1}{2}}, \Delta y_j, \Delta y_{j-1}, \Delta y_{j+1}), \\
		     \bar u^{(+)}_{i,j-\frac12} &= \tilde{\bar u}_{i,j} - \tfrac{1}{2}\; H(\delta_{i,j+\frac{1}{2}},\delta_{i,j-\frac{1}{2}}, \Delta y_j, \Delta y_{j+1}, \Delta y_{j-1}).			     
		    \end{split}
		    \end{align}
	\end{subequations}
  	Here, $H = H_{3, \text{neq}}$ or $H = H^{(c)}_{3L,\text{neq}}$.
  	The two-dimensional undivided differences $\delta_{i\pm\frac{1}{2},j}, \delta_{i, j\pm\frac{1}{2}}$, 
  	are defined as above, Eq.~\eqref{eq:deltas2D} and are adapted to the non-uniform setting as explained 
  	in Sec.~\ref{sec:nonUniform}.
\end{enumerate}	
Steps 2.-5. remain the same, see Sec.~\ref{sec:2DCartesian}.
% %----------------------------------------------------------------------------------------------------------------------------------------------------------------------
\section{Numerical examples}\label{sec:numerics} %----------------------------------------------------------------------------------------------------------------------------------------------------------------------
In this section we present different numerical examples validating the concepts introduced in Sec.~\ref{sec:limiter2D}. We first prove in Sec.~\ref{sec:numerics1DnonUniform} that third-order accuracy is 
obtained on non-equidistant one-dimensional grids. In Sec.~\ref{sec:numerics2Dvortex}, the vortex evolution 
performed on a two-dimensional Cartesian mesh shows that also in 2D, the limiter yields third-order accuracy. 
Then, Sec.~\ref{sec:numerics2Dadvection} presents the two-dimensional advection equation on a non-uniform 
Cartesian mesh. Finally, the double Mach reflection, Sec.~\ref{sec:numericsDMR} and the two-dimensional 
Riemann problem with four shocks show the excellent performance of $H_{3\text{L}}^c$ using AMR.

All simulations were performed using the third-order accurate strong stability preserving Runge-Kutta (SSP-RK3) 
time integrator developed by Gottlieb et.~al.~\cite{GottliebShuTadmor2001}.
\newpage
%
%----------------------------------------------------------------------------------------------------------------------------------------------------------------------
\subsection{Testing the Convergence Order on a Non-Uniform 1D Grid}\label{sec:numerics1DnonUniform}
%----------------------------------------------------------------------------------------------------------------------------------------------------------------------
%
\begin{figure}[t]
	\centering
	\includegraphics[width=\textwidth]{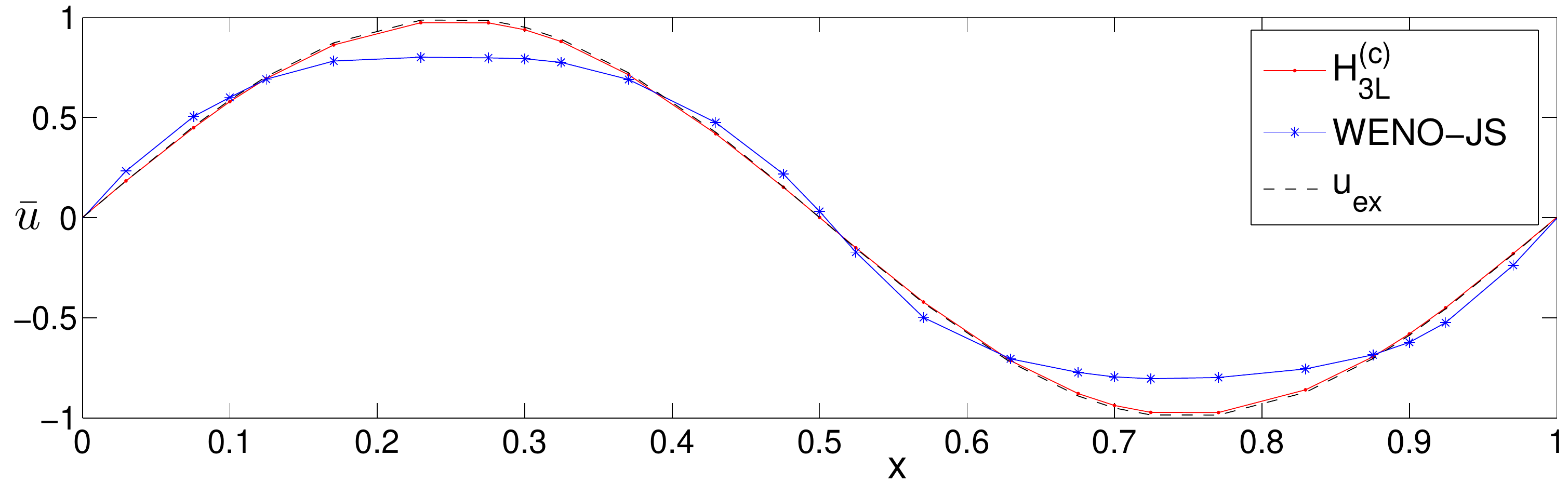}			
	\caption{Comparison of solutions obtained with WENO-JS and 
          $H_{3\text{L}}^c$. Test case \eqref{eq:smoothAdvEq} with $N=25$ grid
          cells on $[0,1]$, CFL $0.95$, until time $t_\text{end}=1$.} 
  	\label{fig:comparisonUnstructured1D}  
\end{figure}
\FloatBarrier
In this section we want to verify that the extension of the third-order limiter 
function $H_{3\text{L}}^c$ from equidistant 
to non-equidistant grids still yields third order accurate solutions. Thus, we 
consider the linear advection equation with smooth initial conditions
\begin{align}
	\begin{cases}
		u_t+ u_x &= 0\\
		u(x,0) &= \sin(2\pi\,x)
	\end{cases}
	\label{eq:smoothAdvEq}
\end{align}
on the domain $[0,1]$ with periodic boundary conditions. In order to verify 
the order of convergence we carry out simulations with $N=25\times 2^j, j=0,\ldots,6$ 
grid cells with end time $t_\text{end}=1.0$ and CFL number $0.95$. Since we are 
interested in non-equidistant grids, the original grid is perturbed by adding 
$c_1\cdot\sin(c_2\,2\pi\,x_{i+1/2})$ to each cell boundary $x_{i+1/2}$ with 
some constants $c_1, c_2\in\mathbb{R}$. In this test case, $c_1 = (10\cdot c_2)^{-1}$ 
and $c_2=5$ have been applied. 

Fig.~\ref{fig:comparisonUnstructured1D} shows the exact solution 
as well as the solution obtained with $H_{3\text{L}}^c$ on a grid with 25 cells. 
This solution is compared to the third order WENO method developed by Liu et.~al.~
\cite{LiuOsherChan1994} with the smootheness measure by Jiang and Shu \cite{JiangShu1996}. 
This scheme is denoted by WENO-JS. The choice of depicting a coarse grid emerges from 
the fact that the non-equidistant mesh structure is well-visible. Also, the improved 
solution quality of the limiter function can be best observed on coarse meshes, as for 
fine grids, all convergent methods look the same.

Finally, to verify that the limiter function is third-order accurate on non-equidistant grids, 
Table \ref{tab:errorsAndEOCofH3Lc} displays the $L_1$- and $L_\infty$-errors of $H_{3\text{L}}^c$. 
The corresponding empirical order of convergence (EOC) is obtained by \\
$\log(err_{j+1}/err_j)/\log(N_j/N_{j+1})$. It can be seen that the limiter function obtains 
the desired accuracy already on coarse meshes.

\begin{table}[t]
	\caption{Errors of $H_{3\text{L}}^c$ in $L_1$- and $L_\infty$-norm and corresponding empirical order of convergence (EOC).}
\subfloat[][Perturbed grid by adding $\frac{1}{50} \sin(10\pi\,x_{i+1/2})$ to each cell boundary $x_{i+1/2}$.]
	{
	\label{tab:errorsAndEOCofH3Lc}	
	\begin{tabular}{ccccc}
		\hline\noalign{\smallskip}
		Grid &$\|u-u_{ex}\|_1$&EOC& $\|u-u_{ex}\|_\infty$&EOC\\ 	
		\noalign{\smallskip}\hline\noalign{\smallskip}
         $25$& 8.322E-03 &       & 1.311E-02 \\ \hline
         $50$&  1.347E-03 & \textbf{2.63} & 2.209E-03&\textbf{2.57} \\ \hline
         $100$&  1.817E-04 & \textbf{2.89} & 2.921E-04&\textbf{2.92} \\ \hline
         $200$&  2.323E-05 & \textbf{2.97} & 3.663E-05&\textbf{3.00} \\ \hline
         $400$&  2.920E-06 & \textbf{2.99} & 4.589E-06&\textbf{3.00} \\ \hline
         $800$&  3.656E-07 & \textbf{3.00} & 5.743E-07&\textbf{3.00} \\ 
         \noalign{\smallskip}\hline
	\end{tabular}	
	}\;
\subfloat[][Non-equidistant grid with random cell boundaries.\\]
	{
	\label{tab:errorsAndEOCofH3Lc2}
	\begin{tabular}{ccccc}
		\hline\noalign{\smallskip}
		Grid &$\|u-u_{ex}\|_1$&EOC& $\|u-u_{ex}\|_\infty$&EOC\\ 	
		\noalign{\smallskip}\hline\noalign{\smallskip}
         $25$  &  6.412E-03 &       			   & 1.004E-02 \\ \hline
         $50$  &  7.888E-04 & \textbf{3.02} & 1.240E-03&\textbf{3.02} \\ \hline
         $100$&  9.844E-05 & \textbf{3.00} & 1.547E-04&\textbf{3.00} \\ \hline
         $200$&  1.234E-05 & \textbf{3.00} & 1.940E-05&\textbf{3.00} \\ \hline
         $400$&  1.535E-06 & \textbf{3.01} & 2.411E-06&\textbf{3.01} \\ \hline
         $800$&  1.930E-07 & \textbf{2.99} & 3.033E-07&\textbf{2.99} \\ 
         \noalign{\smallskip}\hline
	\end{tabular}	
	}
\end{table}
\FloatBarrier
%
%----------------------------------------------------------------------------------------------------------------------------------------------------------------------
\subsection{2D Advection Equation on Non-Uniform Rectangular Grids}\label{sec:numerics2Dadvection}
\begin{figure}[htbp!]  
  \centering
  \subfloat[][Non-uniform mesh with $30\times 30$ grid cells, obtained by 
	Eq.~\eqref{eq:transformationToNonUniform} with $\delta x=\delta y=0.1, c_x=2, c_y=1$.]{
    \includegraphics[width=0.4\textwidth]{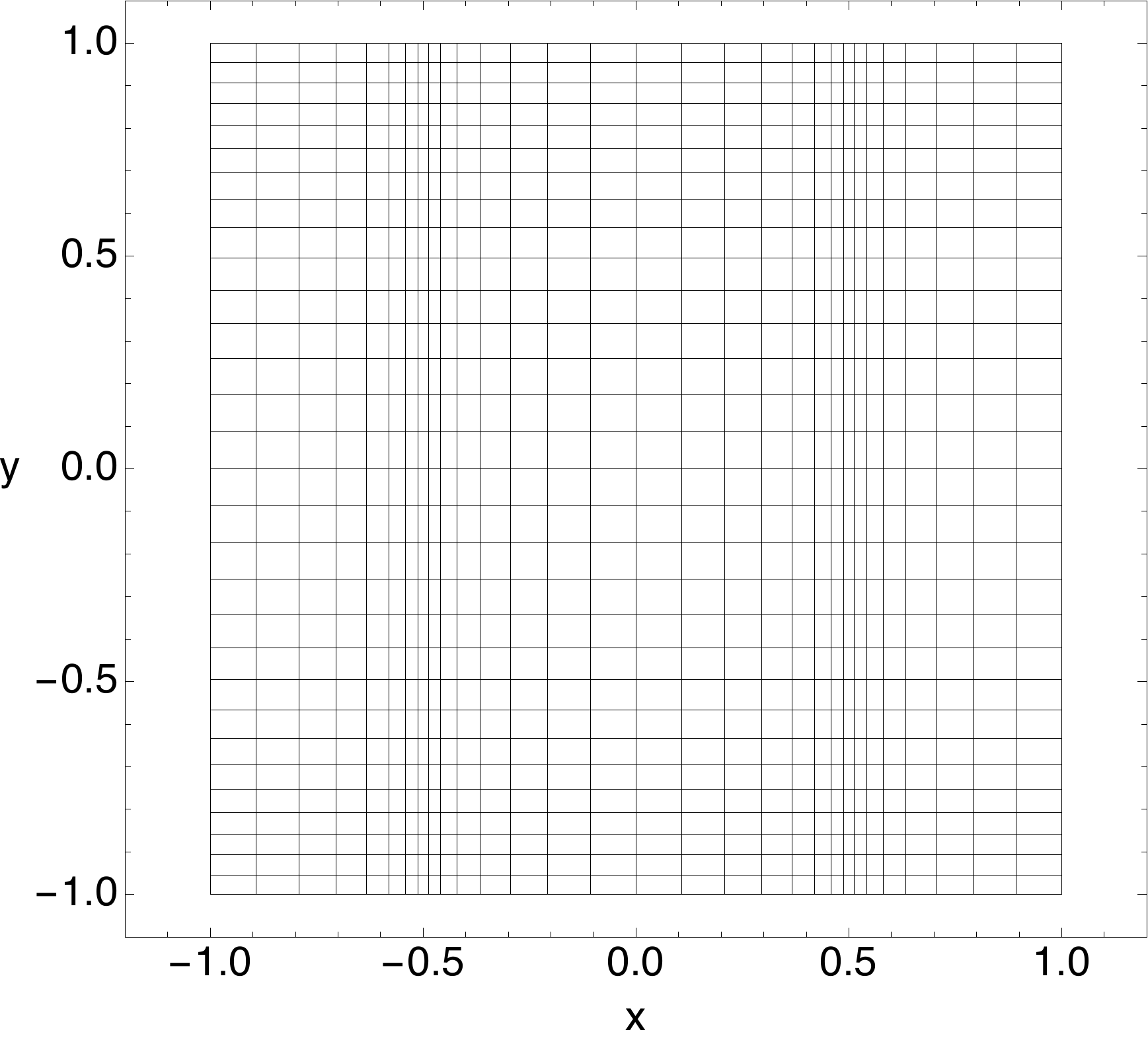}
    \label{fig:nonUniformMesh}
  }\;
  \subfloat[][Solution obtained with \Hlim.]{
	\includegraphics[width=0.54\textwidth]{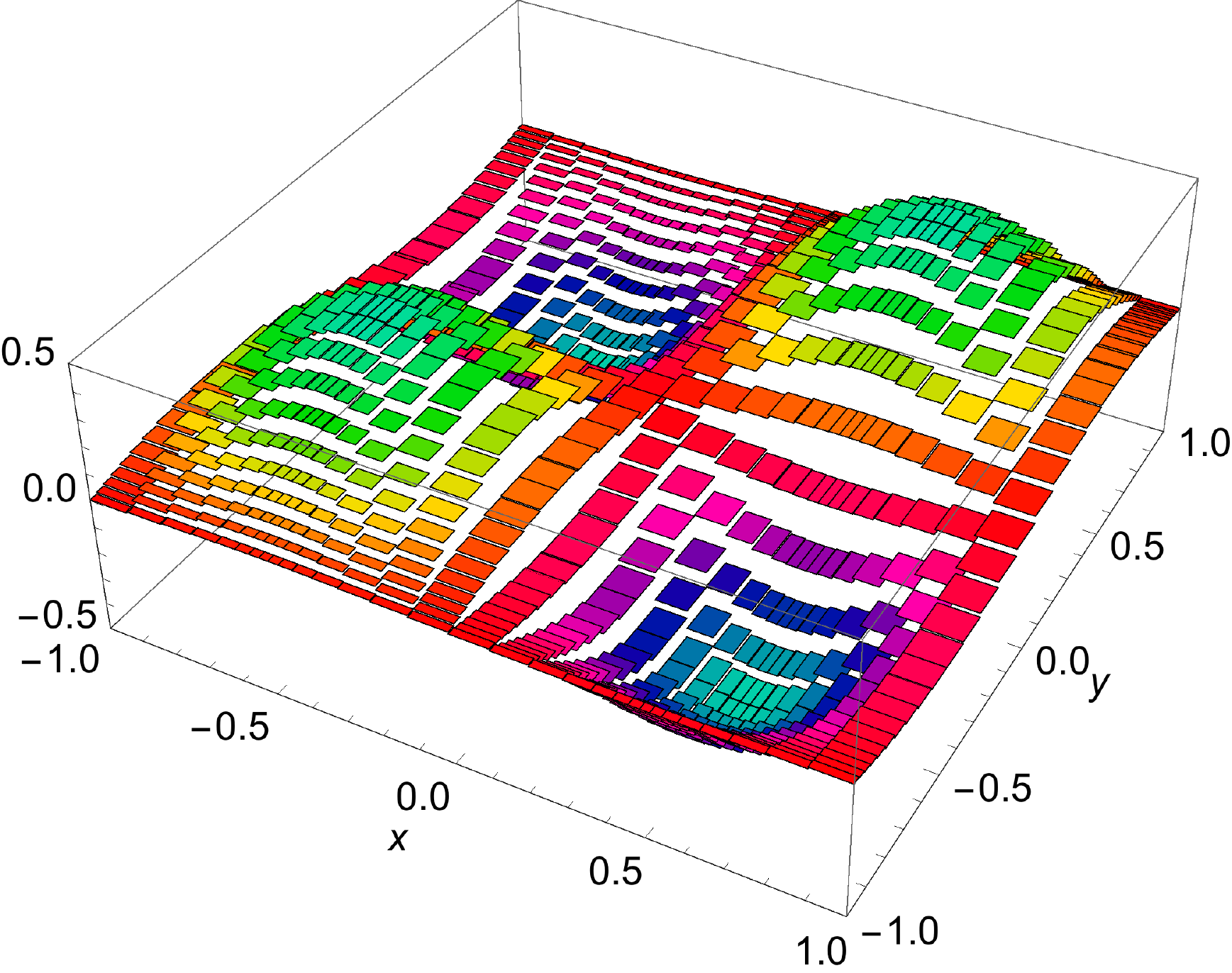}  
	\label{fig:sol30x30uNeq}
  }
  \caption{Computation of Eq.~\eqref{eq:smoothAdvEq2D} on a non-uniform mesh with $30\times 30$ grid cells.}
\end{figure}

%----------------------------------------------------------------------------------------------------------------------------------------------------------------------
%
This numerical problem verifies the accuracy of the two-dimensional numerical scheme on non-uniform 
grids, described in Sec.~\ref{sec:2DnonUniform} and \ref{sec:2DCartesian}. We consider the 
two-dimensional linear advection equation with smooth initial conditions
\begin{align}
	\begin{cases}
		u_t+ a\, u_x + b\,u_y &= 0\\
		u(x, y, 0)=u_0(x, y)  &= \frac12 \sin(\pi\,x)\sin(\pi\,y).
	\end{cases}
	\label{eq:smoothAdvEq2D}
\end{align}
The computational domain is set to $\Omega = [-1, 1]\times [-1, 1]$ and the non-uniformity is obtained 
by Eq.~\eqref{eq:transformationToNonUniform} with $\delta x=0.1, c_x=2, \delta y=0.1, c_y=1$. 
Appying an advection speed of either $(a,b)= (1,0)$ or $(a,b)= (1,1)$ and the simulation time 
$T_{\text{end}} = 2$, the initial condition can be used as exact solution. Thus, the $L_1$-error 
of the numerical solution $u_{ij}^n$ can easily be computed as 
$\|u_{ij}^n - u_0(x_i,y_j)\|_1 = |C_{i,j}|\sum_{i,j}|u_{ij}^n - u_0(x_i,y_j)|$. 
For the simulation, the CFL condition $0.5$ has been imposed and for the decision criterion $\eta$ 
of the limiter function 
$H_{3L,\text{neq}}^{(c)}$, the input value $\alpha = \pi^2$ is obtained by Eq.~\eqref{eq:alpha2D}.
In order to verify the order of convergence, we carry out simulations with 
$N=\{5\times 5, 10\times 10, 20\times 20, 30\times 30, 50\times 50\}$ grid cells. The mesh with 
$30\times 30$ grid cells, perturbed as described above, is depicted in Fig.~\ref{fig:nonUniformMesh} 
and the solution obtained using \Hlim is shown in Fig.~\ref{fig:sol30x30uNeq}. The $L_1$-errors and 
the corresponding empirical orders of convergence are given in Table \ref{tab:errorsAndEOCofH3Lc2D}. The errors for advection speed $(1,0)$ 
are by a mean factor of $0.7$ better than the errors of the solutions advected 
in diagonal direction $(1,1)$. Nevertheless, both simulations yield third order accuracy, see Table \ref{tab:errorsAndEOCofH3Lc2D}
\begin{table}[htb]
	\caption{Errors of $H_{3\text{L}}^c$ in the $L_1$-norm and corresponding empirical order of convergence (EOC).}	
	\label{tab:errorsAndEOCofH3Lc2D}
	\begin{tabular}{ccccc}
	 \hline\noalign{\smallskip}
    &adv. in $x$-dir. &&adv. in diag. dir.\\
		Grid &$\|u-u_{ex}\|_1$&EOC &$\|u-u_{ex}\|_1$&EOC\\ 	
		\noalign{\smallskip}\hline\noalign{\smallskip}
         $5\times 5$    & 3.886E-01 &        		& 5.671E-01&\\ \hline
         $10\times 10$&  1.536E-01 & \textbf{1.34}	& 2.503E-01 & \textbf{1.18}  \\ \hline
         $20\times 20$&  2.621E-02 & \textbf{2.55}  & 6.626E-02 & \textbf{1.92}\\ \hline
         $30\times 30$&  8.304E-03 & \textbf{2.83} 	& 2.088E-02 & \textbf{2.85} \\ \hline
         $50\times 50$&  1.856E-04 & \textbf{2.93} 	& 4.410E-03 & \textbf{3.04} \\
         \noalign{\smallskip}\hline         
	\end{tabular}	
\end{table}
\FloatBarrier
\FloatBarrier
%
%----------------------------------------------------------------------------------------------------------------------------------------------------------------------
\subsection{2D Vortex Evolution}\label{sec:numerics2Dvortex}
%----------------------------------------------------------------------------------------------------------------------------------------------------------------------
%
This problem, originally proposed by Hu \cite{Hu:1999}, describes a two-dimensional vortex
evolution on the periodic domain $[-7,7]\times[-7,7]$, where the flow is described by the
Euler equations. 
The initial data consists of a mean flow $\rho=u=v=p=1$, perturbed by
% adding the small fluctuation
% 
\begin{equation}
  \begin{pmatrix}
    \delta \rho \\\delta u \\\delta v \\\delta p 
  \end{pmatrix}
  =
  \begin{pmatrix}
    (1+\delta T)^{1/(\gamma-1)}-1 \\
    -y\frac{\sigma}{2\pi}e^{0.5(1-r)}\\
    x\frac{\sigma}{2\pi}e^{0.5(1-r)}\\
    (1+\delta T)^{\gamma/(\gamma -1)}-1
  \end{pmatrix}.
\end{equation}
The perturbations in density and pressure are expressed in terms of perturbation in
temperature, $\delta T$, given by
\begin{equation}
  \delta T =-\frac{(\gamma-1)\sigma^2}{8\gamma\pi^2}e^{1-r^2},  
\end{equation}
with $r^2 = x^2 + y^2$, the adiabatic index $\gamma=1.4$, and the vortex strength
$\sigma = 5$. The initial data is also used as a reference solution at time $t=14$,
where it agrees with the exact solution. For the limiter function we need to compute the
radius of the asymptotic region, $\alpha$, given by Eq.~\eqref{eq:alpha2D}. In this case,
$\alpha = 7.9$ is used for the simulation. 
By applying the method on a two-dimensional Cartesian grid, as suggested in 
Sec.~\ref{sec:2DCartesian}, we obtain the full third order, as shown in Table \ref{tab:vortex1}.
The results compare well with the third-order WENO-Z3 implementation developed by Borges 
et.~al.~\cite{Borges:2008} and further improved by Don and Borges \cite{Don:2013}.
Both schemes, the $H^{(c)}_{3L}$ limiter function and WENO-Z3, are implemented 
following the algorithm described in Sec.~\ref{sec:2DCartesian}. 
\begin{table}[t]
  %\centering
  \caption{Results for  the Vortex evolution problem an a uniform grid.}
  \label{tab:vortex1}
  \begin{tabular}[htb]{ccccc}
    \hline\noalign{\smallskip}
    &WENO-Z3&&$H^{(c)}_{3L}$\\
    Grid & $\|\rho-\rho_{ex}\|_1$&EOC& $\|\rho-\rho_{ex}\|_1$&EOC\\
   \noalign{\smallskip}\hline\noalign{\smallskip}
    $64\times 64$      &  1.284E-03 &           & 1.186E-03 & \\ \hline
    $128\times 128$   &   2.302E-04 & \textbf{2.48}& 2.281E-04 & \textbf{2.38} \\ \hline
    $256\times 256$   &   3.125E-05 & \textbf{2.88}& 3.120E-05 & \textbf{2.87} \\ \hline
    $512\times 512$   &   3.955E-06 & \textbf{2.98}& 3.953E-06 & \textbf{2.98} \\ \hline
    $1024\times 1024$ &   4.954E-07 & \textbf{3.00}& 4.954E-07 & \textbf{3.00} \\ \hline
    $2048\times 2048$ &   6.194E-08 & \textbf{3.00}& 6.194E-08 & \textbf{3.00} \\ 
    \noalign{\smallskip}\hline
  \end{tabular}  
\end{table} 
\FloatBarrier
% ----------------------------------------------------------------------------------------------------------------------------------------------------------------------
\subsection{Double Mach Reflection}\label{sec:numericsDMR}
%----------------------------------------------------------------------------------------------------------------------------------------------------------------------
%
In this section we apply the limiter function on a Cartesian grid with AMR, as described in 
Sec.~\ref{sec:AMR}. The test case consists of the double Mach reflection problem 
proposed by Woodward and Colella \cite{Woodward:1984}. It describes a Mach 10 
shock reflection off a 30--degree wedge.
The computational domain is the rectangle $[0,3] \times [0,1]$. 
To obtain the same resolution for both, the $x-$ and $y-$direction, each block contains
$36 \times 12$ mesh cells.  
We set level 3 as the coarsest level and allow up to 4 additional refinements, thus the
finest level   corresponds to a discretization with $4608 \times 1536$ mesh cells. 
The refinement threshold is set to $\delta_0 = 2000$. 
Due to the constant initial date, $\alpha$ turns out to be $0$, cf.~Eq.~\eqref{eq:alpha2D}. 
Therefore, the combined limiter function \Hlim, Eq.~\eqref{eq:H3Lc}, reduces to $H_{3\text{L}}$,
see Eq.~\eqref{eq:limiterH3L}.

\begin{figure}[t!]
  \centering
    \includegraphics[width=\textwidth]{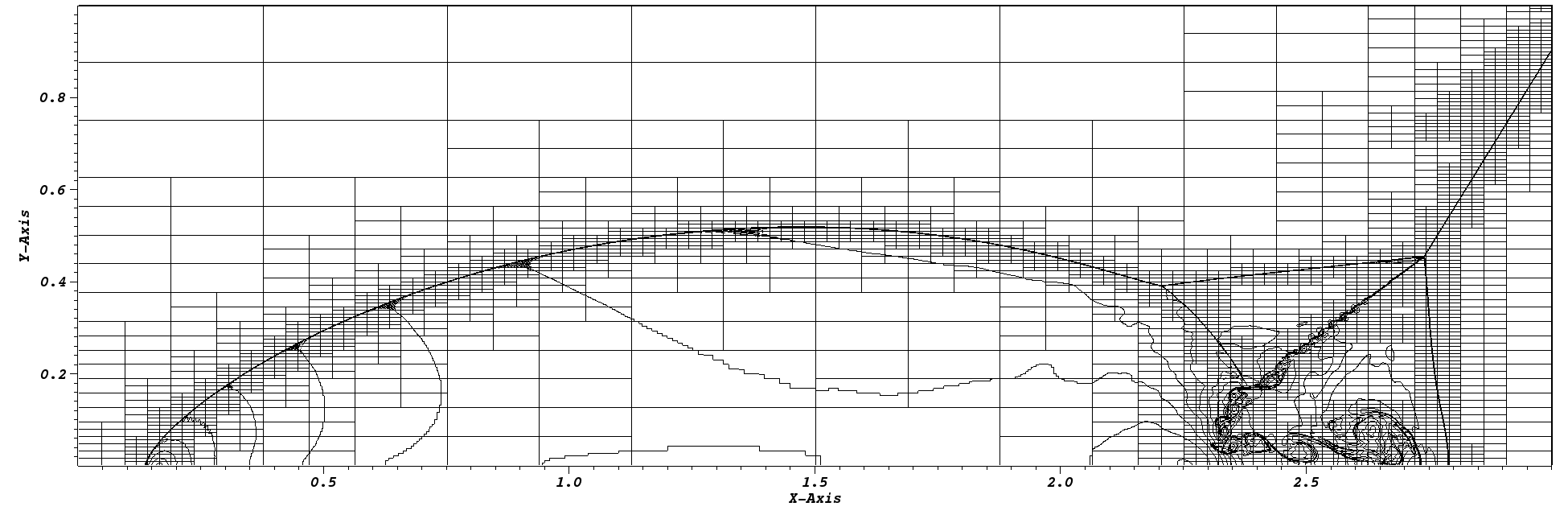}
  \caption{\label{fig:DMR_new} 
    AMR computation of the Double Mach Reflection problem computed with $H^{(c)}_{3L}$. 
    The grid resolution corresponds to a grid with $288\times 96$ mesh cells on the
    coarsest level and up to $4608 \times 1536$ mesh cells on the finest level. 
  } 
\end{figure}

\begin{figure}[htbp!]  
  \centering
  \subfloat[][Computed with WENO-Z3 on an adaptive grid of level 3--7.]{
    \includegraphics[width=0.47\textwidth]{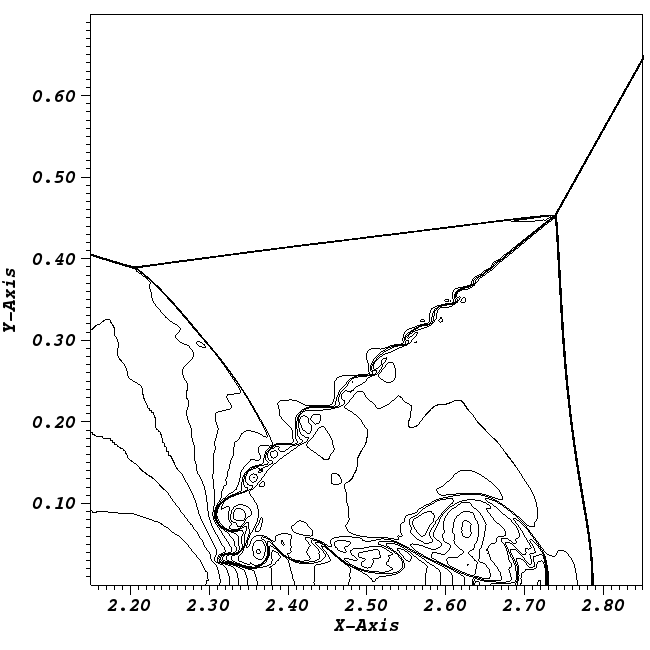}
    \label{fig:DMRweno4}
  }\;
  \subfloat[][Computed with $H^{(c)}_{3L}$ on an adaptive grid of level 3--7.]{%
    \includegraphics[width=0.47\textwidth]{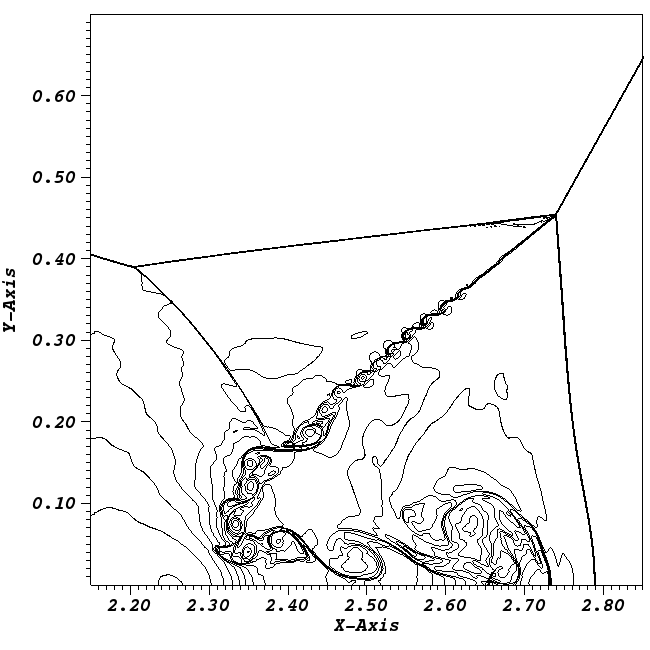}
    \label{fig:DMRh4}
  }\\
  \subfloat[][Computed with WENO-Z3 on an adaptive grid of level 3--8.]{
    \includegraphics[width=0.47\textwidth]{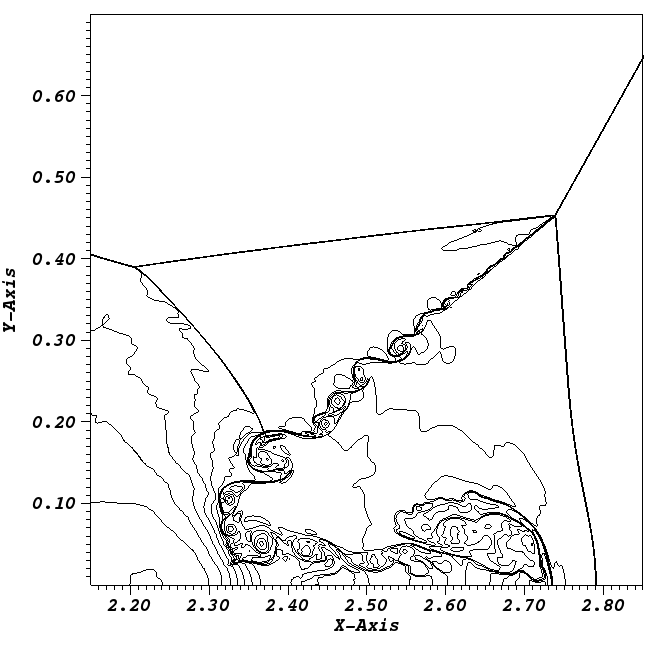}
    \label{fig:DMRweno5}
  }\;
  \subfloat[][Computed with $H^{(c)}_{3L}$ on an adaptive grid of level 3--8.]{%
    \includegraphics[width=0.47\textwidth]{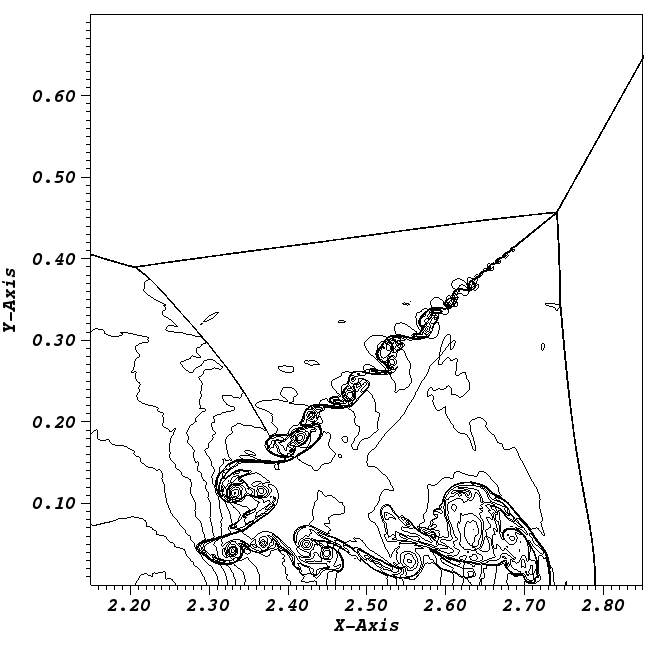}
    \label{fig:DMRh5}
  }
  \caption{\label{fig:DMR_zoom} 
   AMR computation of the Double Mach Reflection problem.}  
\end{figure}

Fig.~\ref{fig:DMR_new} shows the result of the simulation using the third-order limiter 
function $H_{3L}$ at time $t_{end}=0.2$, including the block structure.
\FloatBarrier
A close--up view of the Mach stem region is shown in Fig.~\ref{fig:DMR_zoom}.  
The computations were performed with four (Fig.~\ref{fig:DMRweno4} and \ref{fig:DMRh4}) 
and five (Fig.~\ref{fig:DMRweno5} and \ref{fig:DMRh5}) levels of refinement.
For comparison, we also show the results of a third-order WENO-Z reconstruction on 
the left hand side. 
In direct comparison the $H^{(c)}_{3L}$ scheme produces more roll ups in the inner
region. 
This is a desired feature since the slip line is physically instable, indicating that 
the scheme with \Hlim introduces less numerical viscosity than WENO-Z3.   
\FloatBarrier
%
%-----------------------------------------------------------------------
\subsection{2D Riemann Problem}\label{sec:numerics2dRP}
%-----------------------------------------------------------------------
%
\begin{figure}[t!]  
  \centering
  \subfloat[][Computed with WENO-Z3 on a uniform grid.\\ $1024\times 1024$ mesh cells.]{
    \includegraphics[width=0.46\textwidth]{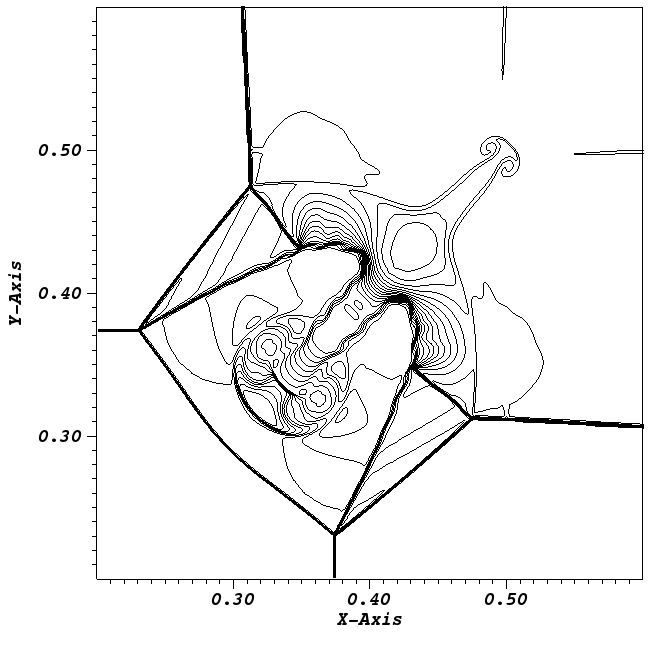}
    \label{fig:RP_weno1024}
  }\; 
  \subfloat[][Computed with $H^{(c)}_{3L}$ on a uniform grid.\\ $1024\times 1024$ mesh cells. ]{%
    \includegraphics[width=0.46\textwidth]{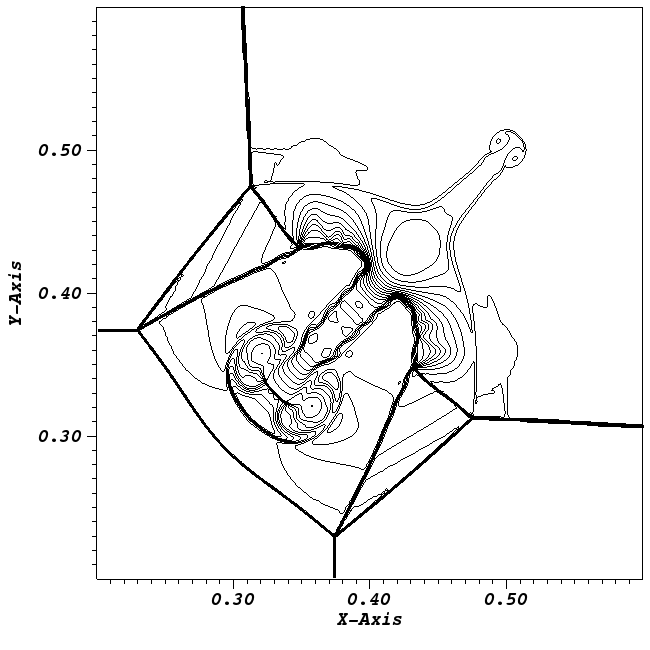}
    \label{fig:RP_h1024}
  }\\
  \subfloat[][Computed with WENO-Z3 on a uniform grid.\\$2048\times 2048$ mesh cells.]{
    \includegraphics[width=0.46\textwidth]{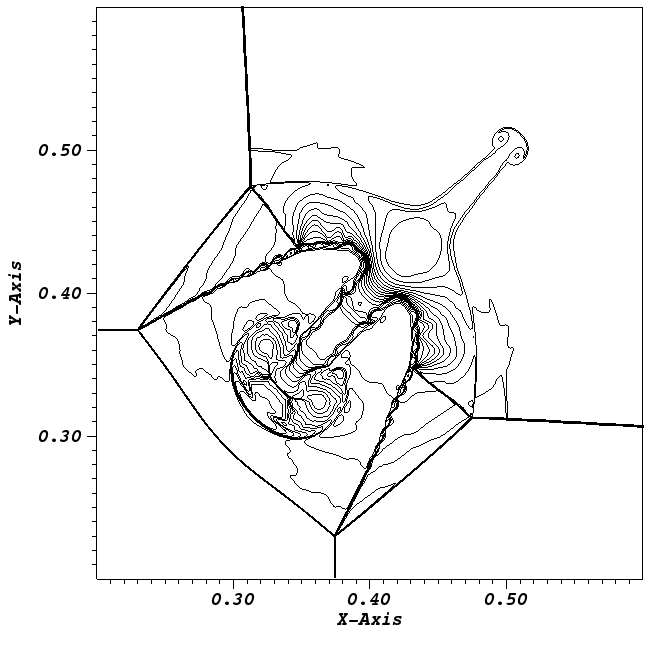}
    \label{fig:RP_weno2048}
  }\;
  \subfloat[][Computed with $H^{(c)}_{3L}$ on a uniform grid.\\ $2048\times 2048$ mesh cells. ]{%
    \includegraphics[width=0.46\textwidth]{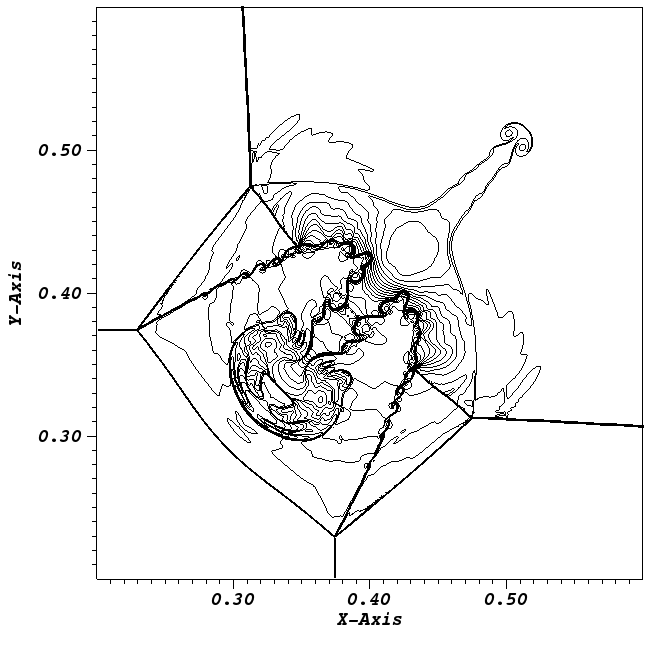}
    \label{fig:RP_h2048}
  }
  \caption{\label{fig:rp} 
   Results of the 2D Rimann problem at final time  $t_{end} = 0.3$.}  
\end{figure}
\FloatBarrier
The next testcase we consider is a configuration of four interacting shocks in the domain $[0,1]^2$. 
The initial values have the form

  \begin{equation}\label{eqn:2drp}
    (\rho, p, u, v)(x,y,0) =
    \left\{ \begin{array}{c c}
        (1.5      , 1.5    , 0.0      , 0.0      ) & x>0.5, y>0.5\\
        (0.5323, 0.3    , 1.2060, 0.0      ) & x<0.5, y>0.5\\
        (0.1380, 0.029, 1.2060, 1.2060) & x<0.5, y<0.5\\
        (0.5323, 0.3    , 0.0      , 1.2060) & x>0.5, y<0.5 
      \end{array}\right..
  \end{equation}  

This testcase was originally proposed by Schulz-Rinne 
\cite{Schulz-Rinne:1993,Schulz-Rinne:1993a} along with several other configurations of 2D
Riemann problems. 
Due to the constant initial data we set here $\alpha  = 0$ for the limiter \Hlim, see \ref{sec:AppendixA}.

Fig.~\ref{fig:rp} shows the results at final time $t_{end} = 0.3$. 
The results obtained by applying the limiter function \Hlim are compared to results computed with 
the third-order WENO-Z reconstruction. As for the double Mach reflection, Sec.~\ref{sec:numericsDMR}, 
the scheme with \Hlim introduces less numerical viscosity than WENO-Z3. 
%
%----------------------------------------------------------------------------------------------------------------------------------------------------------------------
\section{Conclusion}\label{sec:Conclusion}
%----------------------------------------------------------------------------------------------------------------------------------------------------------------------
%
In this work we have extended the recently proposed third-order limiter function $H^{(c)}_{3L}$ 
\cite{SchmidtmannSeiboldTorrilhon2015} from one-dimensional equidistant grids to 
non-uniform and Cartesian AMR meshes in two space dimensions. For the reconstruction of 
interface values, the presented limiter function takes into account the smallest possible stencil 
for reaching third-order accuracy. Thus, in one space dimension the reconstruction 
necessitates three cell mean values and in two dimensions five cells. 

For the limiter to be applicable to one dimensional non-equidistant grids, the undivided differences 
$\delta_{i- 1/2} = \bar{u}_{i}-\bar{u}_{i-1}$ and $\delta_{i+1/2} = \bar{u}_{i+1}-\bar{u}_{i}$ 
have been adapted to be meaningful again. The resulting expressions are closely 
related to the smootheness indicators by Jiang and Shu \cite{JiangShu1996}. A numerical test 
case verifies that the resulting scheme yields the desired third-order accuracy.

For the two dimensional test cases, the popular approach of dimension splitting has been applied. 
In order to use this method without loss of third-order accuracy, we use the order-fix developed in 
\cite{BuchmuellerHelzel2014}. Here, the scheme has first been extended to Cartesian grids, then 
we showed that it can be incorporated in a scheme with adaptive mesh refinement and also 
two-dimensional, non-uniform rectangular grids have been proven to yield third-order accurate 
solutions.

The resulting scheme has been tested on a number of numerical examples and 
shows the desired third-order accuracy. We also compared the limiter function to
the third-order WENO-Z3 reconstruction \cite{Borges:2008, Don:2013} 
and obtain equally good results. \Hlim seems to introduce less viscosity into the scheme, thus 
being able to better reproduce the physically instable details of the double Mach reflection problem.

%The proposed finite volume limiter function significantly reduces oscillations near discontinuities 
%which could now also be verified on more general grids.
%
%==============================
\begin{appendix}
\section{Appendix}
\subsection{Brief Summary of $H_{3\text{L}}^{(c)}$ on Uniform-Grids}\label{sec:AppendixA}
	In this section we want to recall the third-order limiter function developed in 
	\cite{SchmidtmannSeiboldTorrilhon2015}. In one space dimension on equidistant 
	grids it reads
	\begin{align}
		\label{eq:limiterH3L}
		H_{3\text{L}}(\delta_{i-\frac{1}{2}},\delta_{i+\frac{1}{2}}) = 
		\,&\sgn(\delta_{i+\frac{1}{2}})\;\max(0,\min(\sgn(\delta_{i+\frac{1}{2}})\,H_{3}, \max(-\sgn(\delta_{i+\frac{1}{2}})\delta_{i-\frac{1}{2}}, \\
		&\min(2\,\sgn(\delta_{i+\frac{1}{2}})\delta_{i-\frac{1}{2}}, \sgn(\delta_{i+\frac{1}{2}})\,H_{3}, 1.5 |\delta_{i+\frac{1}{2}}|)))). \nonumber
	\end{align}
	As described in \cite{SchmidtmannSeiboldTorrilhon2015}, there exist cases, where the limiter 
	function $H_{3\text{L}}$ is unable to distinguishing between smooth extrema and discontinuities
	due to the constraint of using three cells only. Therefore, the decision criterion $\eta$ 
	was introduced, which is able to distinguish between smooth extrema and discontinuities in most 
	cases. The switch function is defined by
	\begin{subequations}
		\begin{align}
			\label{eq:eta}
			\eta=\eta(\delta_{i-\frac{1}{2}},\delta_{i+\frac{1}{2}}) = 
			\frac{\sqrt{(\delta_{i-\frac12})^2+(\delta_{i+\frac12})^2}}{\sqrt{\frac{5}{2}}\,\alpha\,\Delta x^2}, 
			\intertext{where $\alpha$ denotes the maximum second derivative of the initial conditions}
 			\label{eq:alpha}
 			\alpha \equiv \max_{x \in \Omega \backslash \Omega_d} |u_{0}^{\prime\prime}(x)|.
 		\end{align} 		 		
	\end{subequations}	
	Here, $\Omega$ is the computational domain and $\Omega_d\subset \Omega$ the subset containing 
	discontinuities. This means that the maximum second derivative is only considered in smooth 
	parts of the domain. With the switch function $\eta$, the combined limiter function reads
	\begin{align}
        \label{eq:check}
		H_{3\text{L}}^{(c)} (\delta_{i-\frac{1}{2}},\delta_{i+\frac{1}{2}}) = 
			\begin{cases}
				H_{3}(\delta_{i-\frac{1}{2}},\delta_{i+\frac{1}{2}}) \quad\qquad &\text{if}\;\eta <1  \\
		  	  	H_{3\text{L}}(\delta_{i-\frac{1}{2}},\delta_{i+\frac{1}{2}})\quad\; &\text{if}\;\eta \geq 1.
		  	\end{cases}
	\end{align}	
	Note that for performance reasons, instead of computing $\eta$ in every cell by using \eqref{eq:check} we precompute 
	\begin{equation*}
	  \tau = \frac{5}{2} \left(\alpha\Delta x^2\right)^2 
	\end{equation*}
	and use
	\begin{align}
	  \label{eq:check_racoon}
	  H_{3\text{L}}^{(c)} (\delta_{i-\frac{1}{2}},\delta_{i+\frac{1}{2}}) = 
	  \begin{cases}
	    H_{3}(\delta_{i-\frac{1}{2}},\delta_{i+\frac{1}{2}}) \quad 
	    &\text{if}\;\delta^2_{i-\frac12}+\delta^2_{i+\frac12}<\tau  \\
	    H_{3\text{L}}(\delta_{i-\frac{1}{2}},\delta_{i+\frac{1}{2}})\quad\; 
	    &\text{if}\;\delta^2_{i-\frac12}+\delta^2_{i+\frac12}\geq\tau
	  \end{cases}
	\end{align}
	instead of Eq.~\eqref{eq:check}. This leads to a significant reduction of computational time. 
\subsection{Derivation of the Full-Third-Order Reconstruction on Non-Equidistant Grids}\label{app:full3rdOrderNeq}
	The quadratic polynomial $p_i(x)$ which satisfies
	\begin{align*}
		\frac{1}{\Delta x_{i+\ell}}\int_{x_{i-\frac12+\ell}}^{x_{i+\frac12+\ell}} p_i(x) dx = \bar{u}_{i+\ell},\quad \ell\in \{-1, 0, 1\}
	\end{align*}
	is given by
	\begin{subequations}
	\begin{align}
	p_i(x) &= a\,(x-x_i)^2 + b\,(x-x_i) + c \\
	\intertext{with}
	a &= \frac{1}{2}\frac{ 
		\frac{\Delta x_i+\Delta x_{i+1}}{2}(u_{i-1}-u_i) + \frac{\Delta x_i+\Delta x_{i-1}}{2}(u_{i+1}-u_i) 
		}{ 
		\left(\frac{\Delta x_i+\Delta x_{i-1}}{2}\right)\left(\frac{\Delta x_i+\Delta x_{i+1}}{2}\right)\Delta_i 
		},\\
	%--------------------------------------------------------
	b &= \frac{
	u_i (\Delta x_{i+1}-\Delta x_{i-1}) (2 \Delta_i+3\Delta x_i) 
	}{ 
		\left(\frac{\Delta x_i+\Delta x_{i-1}}{2}\right)\left(\frac{\Delta x_i+\Delta x_{i+1}}{2}\right)\Delta_i
	} \nonumber \\
	& \quad - 
	\frac{
	\frac{2}{3} u_{i-1} \left(\frac{\Delta x_i+\Delta x_{i+1}}{2}\right) \left(2\frac{(\Delta x_i+\Delta x_{i+1})}{2}+\Delta x_{i+1}\right)
	}{ 
		\left(\frac{\Delta x_i+\Delta x_{i-1}}{2}\right)\left(\frac{\Delta x_i+\Delta x_{i+1}}{2}\right)\Delta_i
	} \nonumber \\
	&\quad + \frac{
	\frac{2}{3} u_{i+1} \left(\frac{\Delta x_i+\Delta x_{i-1}}{2} \right) \left(2\frac{(\Delta x_i+\Delta x_{i-1})}{2}+\Delta x_{i-1}\right)
	}{ 
		\left(\frac{\Delta x_i+\Delta x_{i-1}}{2}\right)\left(\frac{\Delta x_i+\Delta x_{i+1}}{2}\right)\Delta_i
	},\\
	%--------------------------------------------------------
	c &= \frac{
	u_i \left(\Delta x_i\left(6\Delta x_i^2+9\Delta x_i (\Delta x_{i-1}+\Delta x_{i+1})+4 (\Delta x_{i-1}+\Delta x_{i+1})^2\right)+4\Delta x_i\Delta x_{i-1}\Delta x_{i+1} \right)
	}{
	4 (\Delta x_i+\Delta x_{i-1}) (\Delta x_i+\Delta x_{i+1}) 3\Delta_i
	} \nonumber \\
	& \quad+
	\frac{4\Delta x_{i-1}\Delta x_{i+1} (\Delta x_{i-1}+\Delta x_{i+1})u_i - 
     \Delta x_i^2 (u_{i-1} (\Delta x_i+\Delta x_{i+1})+u_{i+1} (\Delta x_i+\Delta x_{i-1}))
	}{
	4 (\Delta x_i+\Delta x_{i-1}) (\Delta x_i+\Delta x_{i+1}) 3\Delta_i
	}.
	\end{align}
\end{subequations}	
Evaluating this polynomial at $x_{i+\halb}$ and rearranging yields the formulation 
\begin{align*}
	u^{(-)}_{i+\frac12}= \bar u_i + \tfrac{1}{2}\; H_{3,\text{neq}}\left(\delta _{i-\frac{1}{2}},\delta _{i+\frac{1}{2}},	\Delta x_i, \Delta x_{i-1}, \Delta x_{i+1}\right)
\end{align*}
with $H_{3,\text{neq}}$ given by 
\begin{align*}
	H_{3,\text{neq}}\left(\delta _{i-\frac{1}{2}},\delta _{i+\frac{1}{2}},	
	\Delta x_i, \Delta x_{i-1}, \Delta x_{i+1}\right) &= \nonumber \\
	\frac{\Delta x_i}{\Delta_i} \frac{1}{3}&\left(2 \frac{\Delta_{i-\frac12}}
	{\Delta_{i+\frac12}}\delta _{i+\frac{1}{2}} + \frac{\Delta x_{i+1}}
	{\Delta_{i-\frac12}}\delta _{i-\frac{1}{2}}\right)
\end{align*}
with the abbreviations
\begin{align*}
	\Delta_i = \frac{\Delta x_{i-1}+\Delta x_i+\Delta x_{i+1}}{3}, \;\,\Delta_{i-\frac12} = 
	\frac{\Delta x_{i-1}+\Delta x_{i}}{2},\;\,\Delta_{i+\frac12} = \frac{\Delta x_i+\Delta x_{i+1}}{2}.
\end{align*}
Evaluating and rearranging $p_i(x_{i-\halb})$ yields
\begin{align*}
	u^{(+)}_{i-\frac12}= \bar u_i - \tfrac{1}{2}\; H_{3,\text{neq}}\left(\delta _{i+\frac{1}{2}},\delta _{i-\frac{1}{2}},	\Delta x_i, \Delta x_{i+1}, \Delta x_{i-1}\right)
\end{align*}
with the same reconstruction function $H_{3,\text{neq}}$. Note however, that the order of the arguments 
has changed in this case leading to
\begin{align*}
	H_{3,\text{neq}}\left(\delta _{i+\frac{1}{2}},\delta _{i-\frac{1}{2}},	
	\Delta x_i, \Delta x_{i+1}, \Delta x_{i-1}\right) &= \nonumber \\
	\frac{\Delta x_i}{\Delta_i} \frac{1}{3}&\left(2 \frac{\Delta_{i+\frac12}}
	{\Delta_{i-\frac12}}\delta_{i-\frac{1}{2}} + \frac{\Delta x_{i-1}}
	{\Delta_{i+\frac12}}\delta_{i+\frac{1}{2}}\right).
\end{align*}
\end{appendix}
%==============================

%----------------------------------------------------------------------------------------------------------------------------------------------------------------------%
%----------------------------------------------------------------------------------------------------------------------------------------------------------------------%
\section*{References}
\bibliographystyle{apa}     
\bibliography{./references_complete}

\end{document}